# Seminormal representations of Weyl groups and Iwahori-Hecke algebras

**Arun Ram***

Department of Mathematics
University of Sydney
NSW 2006 Australia

## ABSTRACT

The purpose of this paper is to describe a general procedure for computing analogues of Young's seminormal representations of the symmetric groups. The method is to generalize the Jucys-Murphy elements in the group algebras of the symmetric groups to arbitrary Weyl groups and Iwahori-Hecke algebras. The combinatorics of these elements allows one to compute irreducible representations explicitly and often very easily. In this paper we do these computations for Weyl groups and Iwahori-Hecke algebras of types $A_n$, $B_n$, $D_n$, $G_2$. Although these computations are in reach for types $F_4$, $E_6$, and $E_7$, we shall, in view of the length of the current paper, postpone this to another work.

## 0. INTRODUCTION

This paper is an attempt to simplify, unify, review, and expand the theory of seminormal representations of the symmetric groups. It is my hope that the material in this paper will be clear enough and general enough to dispel some of the mystique of this subject. I would suggest that the reader begin with section 3, the symmetric group case, and go back and pick up the generalities from sections 1 and 2 as they are needed. This will make

* Research supported in part by a National Science Foundation grant DMS-9300523 at the University of Wisconsin.

1991 *Mathematics Subject Classification.* Primary 20F55, 20C15 Secondary 20C30, 20G05.

the motivation for the material in the earlier sections much more clear and the further examples in the later sections very easy.

In such a subject, where the credit for the discoveries has been traditionally complicated and mixed up, it is difficult to claim credit for much of anything. I believe that my general approach to seminormal representations of Weyl groups and Iwahori-Hecke algebras given in sections 1 and 2 of this paper is new, at least nobody explained it to me before. If this was the structure motivating the work of Hoefsmit [H], then he certainly didn't give any hint of it, he produced the seminormal representations out of thin air and then proved that they were correct. I believe that my realization that the Jucys-Murphy elements are coming from the very natural central elements in (2.1) and Proposition (2.4) is new. I also think that before this paper nobody has realized that there is a simple concrete connection (Proposition (2.8)) between Jucys-Murphy type elements in Iwahori-Hecke algebras and Jucys-Murphy elements in group algebras of Weyl groups (even for the special cases). I know that the analogues of the Jucys-Murphy elements in Weyl groups of types B and D will be new to some of the experts and known to others. These Jucys-Murphy elements for types B and D are not new, similar elements appear in the paper of Cherednik [Ch], but I was not able to recognize them there until they were pointed out to me by M.L. Nazarov. I extend my thanks to him for this. Some people were asking me for Jucys-Murphy elements in type $G_2$ as late as June of 1995. In July 1995 I was told that it was unknown how to quantize the elements of Cherednik, i.e. find analogues of them in the Iwahori-Hecke algebras of types B and D. Of course, this has been done already in 1974, by Hoefsmit.

The idea behind my approach to seminormal representations is as old as representation theory, it is the same method that gives us Gelfand-Tsetlin bases, and that is often used to find the representations of centralizer algebras and operator algebras.

I would also like to point out that there is a theory of seminormal representations and Jucys-Murphy elements for centralizer algebras as well. In the case of the Brauer centralizer algebra the analogues of the Jucys-Murphy elements are due to M.L. Nazarov [Na], and in the case of the Birman-Wenzl algebra to N. Reshetikhin [Re] and R. Leduc and the author [LR]. In the theory of centralizer algebras the Jucys-Murphy elements "come from" the Casimir element of the centralizing Lie algebra or quantum group.

*The mystique*

Unfortunately, the theory of seminormal representations of the symmetric group is a theory which, over the years, has developed quite a bit of mystique. I suppose that this began with the fact that Alfred Young's original work was hard to read and because he was somewhat "in the wrong place at the right time" in that his work was competing with the work of Frobenius. Today there are very few people that have carefully read Young's work and even those people complain about the writing, the notation, and the general mysteriousness. The second major complicating factor was that Hoefsmit's thesis was never published. The result of this was that several people rediscovered parts of his work and giving credit for who did what when became very difficult.

Jucys, and later Murphy, introduced some very special elements that had the chance of making everything in theory of seminormal representations much more simple and clear. This was not to be. The original work of Jucys ([Ju1-3] 1966, 1971, 1974) was published mostly in Lithuanian physics journals and was not read by many in the western mathematical community. Hoefsmit ([H] 1974), unknowingly, had analogues of these elements in the Iwahori-Hecke algebras of type B but his thesis was never published and these elements remained largely unknown. Only, when the work of Murphy ([M1-2] 1981, 1983) appeared did these elements begin to receive wider attention. Somehow though, the mechanism of seminormal representations and Jucys-Murphy elements still seemed complex and specific to the symmetric group.

After this, people began to try to generalize this picture to other types. Cherednik [Ch] gave analogues of the Jucys-Murphy elements for the Weyl groups of types B and D, but his paper was read by almost noone in the seminormal representations camp since his paper was written from the point of view of constructing monodromy representations. In the period 1985-1995, Dipper, James, and Murphy ([DJ1-3] 1986, 1987, 1992, [M1],[M4] 1992, 1995, [DJM] 1995)) have done a lot of work on representations of Iwahori-Hecke algebras and have produced analogues of the Jucys-Murphy elements for Iwahori-Hecke algebras of types A and B. Their version of these elements in type A was written in such a way that it was not clear to anybody that they were the same as the elements that were in Hoefsmit's thesis!

*History*

Alfred Young [Y] wrote down several ways of producing irreducible representations and minimal idempotents in the group algebra of the symmetric group. This material seems to have been challenging to read from the beginning and it is my impression that today there only a few experts on what Young actually wrote. More recent updates of the theory of Young are given in the books by Rutherford [Ru] and James and Kerber [JK].

Although I have not read Young's work myself, it is generally accepted that Young extended most parts of his work on the symmetric group to the cases of the Weyl groups of types B and D. Other treatments of the representation theory of these groups is in the literature but is sparse and scattered.

In the late 1960's, Jucys [Ju], and later in the early 1980's, Murphy [Mu] inserted a new and beautiful feature in this theory by writing down elements in the group algebra of the symmetric group, which, in Young's irreducible seminormal representations of the symmetric group, are always diagonal matrices. Even better, the diagonal entries of these matrices have an easy combinatorial description. They showed that Young's seminormal representations could be reconstructed from the knowledge of these special elements. Analogues of the Jucys-Murphy elements for Weyl groups of types B and D seem to be first in the work of Cherednik [Ch].

In 1974, Hoefsmit [H] completed a fantastic Ph.D. thesis in which he wrote down irreducible seminormal representations for the Iwahori-Hecke algebras of types A, B, and D. In fact, Hoefsmit had, at least for the type B case, analogues of the Jucys-Murphy elements for Iwahori-Hecke algebras.

More recently, there have been new "Hecke algebras" which have been discovered by Ariki [Ar], Ariki-Koike [AK], and Broué-Malle [BM], which are similar to the Iwahori-Hecke algebras of types B and D. Ariki-Koike [AK] and Ariki [Ar] have shown that Hoefsmit's constructions can be extended to give seminormal representations of these algebras as well. Ariki, Koike, and Broué and Malle have given essentially the same analogues of the Jucys-Murphy elements as Hoefsmit had.

### Acknowledgements

I would like to thank especially Curtis Greene and Persi Diaconis for asking me questions, many helpful conversations and their encouragement of this work. R. Stanley and P. Diaconis tried to get me interested in analogues of Jucys-Murphy elements in Iwahori-Hecke algebras in 1992 and 1993. They had asked me at that time whether I knew analogues of these in any other cases outside of Iwahori-Hecke algebras of Type A and the symmetric group. I only really got excited about Jucys-Murphy elements after Curtis Greene explained the story to me in summer 1994. At that time it still seemed to us that analogues of these elements were unknown for Weyl groups of types B and D.

I would like to thank Rob Leduc for pushing and wanting to find explicit representations of the Brauer and Birman-Wenzl algebras in early 1994. It was while we were doing this work [LR] that we first saw analogues of Jucys-Murphy elements coming out of quantum groups. I would like to thank

G. Benkart, L. Solomon, and A.N. Kirillov for their interest and for many conversations during which many parts of this work and this subject solidified for me. I would like to thank P. Orlik and H. Terao for introducing and explaining the work of Brieskorn-Saito and Deligne on braid groups.

I would like to thank National Science Foundation for continuing support of my work, first with a posdoctoral fellowship and then under the grant DMS-9300523 at the University of Wisconsin. I would also like thank G. Lehrer for many interesting discussions and his hospitality and the Australian Research Council for support under a fellowship at University of Sydney where this paper was written.

## 1. WHAT IS A SEMINORMAL REPRESENTATION?

For convenience and simplicity we shall work over the field $\mathbb{C}$ of complex numbers. Let $\{1\} = G_0 \subseteq G_1 \subseteq \cdots \subseteq G_n = G$ be a chain of finite groups. Let $V^\lambda$ be an irreducible module for $G$. Upon restriction to the group $G_{n-1}$ the module $V^\lambda$ decomposes as a direct sum

$$V^\lambda \cong V^{\mu_1} \oplus \cdots \oplus V^{\mu_k}$$

of irreducible modules for $G_{n-1}$. Similarly each of these summands decomposes into irreducible on restriction to $G_{n-2}$, and so on.

A *seminormal basis* of $V^\lambda$ is a basis $B^\lambda = \{v_L\}$ of $V^\lambda$ that explicitly realizes these decompositions, i.e. there is a partition of $B^\lambda$ into subsets $B^{\mu_1}, \ldots, B^{\mu_k}$ such that if $V^{\mu_i} = \mathbb{C}\text{-}span(B^{\mu_i})$ then

$$V^\lambda = V^{\mu_1} \oplus \cdots \oplus V^{\mu_k}$$

as $G_{n-1}$ modules (note that here there is an = sign rather than only $\cong$). Further, we require that each of the subsets $B^{\mu_i}$ partitions into subsets which realize the decomposition upon restricting to $G_{n-2}$, and so on all the way down the chain. Thus, to specify a seminormal basis one must give, not only the basis of $V^\lambda$ but also the series of partitions. The resulting representation

$$\rho^\lambda \colon G \to M_{d_\lambda}(\mathbb{C}), \qquad d_\lambda = \dim(V^\lambda),$$

of $G$ which is specified by $V^\lambda$ and the basis $B^\lambda$ is a *seminormal representation of $G$ with respect to the chain $G_0 \subseteq \cdots \subseteq G_n = G$.*

The concepts of seminormal bases and seminormal representations apply equally well to any chain of split semisimple algebras $\mathbb{C} \cong H_0 \subseteq H_1 \subseteq \cdots H_n = H$.

*The graph $\Gamma$*

Let $\hat{G}_i$ be an index set for the irreducible representations of $G_i$. Define nonnegative integers $c_\mu^\lambda$, $\mu \in \hat{G}_{i-1}$, $\lambda \in \hat{G}_i$, by the restriction rule from $G_i$ to $G_{i-1}$,

$$V^\lambda \!\downarrow_{G_{i-1}}^{G_i} \cong \bigoplus_{\mu \in \hat{G}_{i-1}} c_\mu^\lambda V^\mu,$$

In other words, upon restriction from $G_i$ to $G_{i-1}$ the irreducible module $V^\mu$, $\mu \in \hat{G}_{i-1}$, appears in the irreducible $G_i$-module $V^\lambda$, $\lambda \in \hat{G}_i$, with multiplicity $c_\mu^\lambda$.

Define a graph $\Gamma$ with

> vertices labeled by the elements of the sets $\hat{G}_i$, and such that
> $\mu \in \hat{G}_{i-1}$ and $\lambda \in \hat{G}_i$ are connected by $c_\mu^\lambda$ edges. (1.1)

The graph $\Gamma$ encodes the restriction rules for the chain $\{1\} = G_0 \subseteq \cdots \subseteq G_n = G$. We shall assume that the unique element in $\hat{G}_0$ is denoted $\emptyset$.

Let $\mu \in \hat{G}_r$ and $\lambda \in \hat{G}_s$ where $r < s$. A *path* from $\mu$ to $\lambda$ is a sequence of $s - r$ edges connecting $\mu$ to $\lambda$,

$$L = (\mu = \lambda^{(r)} \xrightarrow{e_r} \lambda^{(r+1)} \xrightarrow{e_{r+1}} \cdots \xrightarrow{e_{s-1}} \lambda^{(s)} = \lambda), \ \text{ such that } \ \lambda^{(i)} \in \hat{G}_i, \ \ r \leq i \leq s.$$

We distinguish paths which "travel" from $\lambda^{(i)}$ to $\lambda^{(i+1)}$ along different edges. Make the following notations:

$\mathcal{L}(\lambda \to \mu)$ is the set of paths from $\lambda$ to $\mu$;

$\mathcal{L}(\lambda)$ is the set of paths from $\emptyset$ to $\lambda$;

$\mathcal{L}(\lambda \to s)$ is the set of paths from $\lambda$ to any element $\mu \in \hat{G}_s$;

$\mathcal{L}(m)$ is the set of paths from $\emptyset$ to any element $\lambda \in \hat{G}_m$;

$\mathcal{L} = \mathcal{L}(n)$, where $n$ is the total number of groups in the chain $G_0 \subseteq \cdots \subseteq G_n$;

$\Omega(\lambda)$ is the set of pairs $(S, T)$ of paths such that $S, T \in \mathcal{L}(\lambda)$;

$\Omega(m)$ is the set of pairs $(S, T)$ of paths such that $S, T \in \mathcal{L}(\lambda)$ for some $\lambda \in \hat{G}_m$.

In general, by "a path in $\Gamma$" we shall mean an element $L = (\lambda^{(0)} \to \cdots \to \lambda^{(n)}) \in \mathcal{L}$.

*Path algebras*

For each $0 \leq m \leq n$ define a *path algebra* $P_m$ over $\mathbb{C}$ with basis $E_{ST}, (S, T) \in \Omega(m)$ and multiplication given by

$$E_{ST} E_{PQ} = \delta_{TP} E_{SQ}. \tag{1.2}$$

Note that $P_0 \simeq \mathbb{C}$. Each of the algebras $P_m$ is isomorphic to a direct sum of matrix algebras

$$P_m \simeq \bigoplus_{\lambda \in \hat{P}_m} M_{d_\lambda}(\mathbb{C}),$$

where $M_d(\mathbb{C})$ denotes the algebra of $d \times d$ matrices with entries from $\mathbb{C}$ and $d_\lambda = Card(\mathcal{L}(\lambda))$. For each $\lambda \in \hat{G}_m$ define an $P_m$-module by defining

$$V^\lambda = \mathbb{C}\text{-}span\{v_L \mid L \in \mathcal{L}(\lambda)\} \quad \text{and} \quad E_{ST}v_L = \delta_{TL}v_S, \qquad (1.3)$$

for all paths $S, T, L \in \mathcal{L}(\lambda)$. The $A_m$-modules $V^\lambda$, $\lambda \in \hat{G}_m$, realize all of the irreducible $P_m$-modules.

Given a path $T = (\lambda \to \cdots \to \mu)$ from $\lambda$ to $\mu$ and a path $S = (\mu \to \cdots \to \nu)$ from $\mu$ to $\nu$ define

$$T * S = (\lambda \to \cdots \to \mu \to \cdots \to \nu)$$

to be the concatenation of the two paths. Let $r < s$ and for each $\lambda \in \hat{G}_r$ and each pair $(P, Q) \in \Omega(\lambda)$ view the element $E_{PQ} \in P_r$ as an element of $P_s$ by the formula

$$E_{PQ} = \sum_{T \in \mathcal{L}(\lambda \to s)} E_{P*T, Q*T}. \qquad (1.4)$$

This defines, in particular, an inclusion of $P_{m-1}$ into $P_m$ for every $m > 0$. Let $\lambda \in \hat{G}_m$ and let $V^\lambda$ be the irreducible representation of $P_m$ corresponding to $\lambda$ as given in (1.3). Then the restriction of $V^\lambda$ to $P_{m-1}$ decomposes as

$$V^\lambda \Big\downarrow_{P_{m-1}}^{P_m} \simeq \bigoplus_{\mu \overset{e}{\longrightarrow} \lambda} V^\mu,$$

where sum is over all edges $\mu \overset{e}{\longrightarrow} \lambda$ that connect an element $\mu \in \hat{G}_{m-1}$ to the $\lambda \in \hat{G}_m$. The basis vectors $v_L$ form a seminormal basis of the $P_m$-module $V^\lambda$.

*Constructing seminormal representations of $\mathbb{C}G$*

As above, let

$$\{1\} = G_0 \subseteq \cdots \subseteq G_n = G \qquad (1.5)$$

be a chain of finite groups, let $\Gamma$ be the graph which describes the restrictions rules for the inclusions in (1.5) and let $P_m$, $0 \le m \le n$ denote the corresponding path algebras. By construction, the path algebras $P_n$ have natural seminormal representations (1.3) with respect to the inclusions $P_0 \subseteq P_1 \subseteq \cdots \subseteq P_n$. Thus, we should try to find an isomorphism

$$\Phi: \quad \begin{array}{ccc} P_n & \simeq & \mathbb{C}G \\ E_{ML} & \mapsto & e_{ML} \end{array} \qquad \text{such that} \quad \Phi(P_i) = \mathbb{C}G_i, \qquad (1.6)$$

for all $0 \le i \le n$. Given such an isomorphism, irreducible seminormal representations are given by the modules $V^\lambda$ in (1.3) where the action of an element $g \in \mathbb{C}G$ is given by

$$gv_L = \Phi(g)v_L, \qquad (1.7)$$

for all $g \in G$.

Suppose that, for each $1 \le k \le n$,

$$Z_k = \{z_{k,j}\}_{1 \le j \le r_k} \qquad (1.8)$$

is a set of central elements in the group algebra $\mathbb{C}G_k$.

**Lemma 1.9.** *Let $z_{k,j}$ be a central element in $\mathbb{C}G_k$. Let $L = (\lambda^{(0)} \to \cdots \to \lambda^{(n)}) \in \mathcal{L}(n)$ be a path in the graph $\Gamma$ and let $\chi^{\lambda^{(k)}}$ be the irreducible character of $G_k$ indexed by the element $\lambda^{(k)} \in \hat{G}_k$. For any choice of isomorphism $\Phi$ between the path algebra $P_n$ and $\mathbb{C}G$ as in (1.6),*

$$z_{k,j}v_L = c_{k,j}(\lambda^{(k)})v_L, \qquad \text{where} \qquad c_{k,j}(\lambda^{(k)}) = \frac{\chi^{\lambda^{(k)}}(z_{k,j})}{\chi^{\lambda^{(k)}}(1)}.$$

*Proof.* By Schur's lemma any central element $z_{k,j} \in \mathbb{C}G_k$ must act by a scalar multiple of the identity in every irreducible representation of $G_k$. Specifically, $z_{k,j}$ acts by the scalar $c_{k,j}(\mu)$ in the irreducible $G_k$-module indexed by $\mu$. Each of the basis vectors $v_L$, $L = (\lambda^{(0)} \to \cdots \to \lambda^{(n)})$ is in an irreducible $P_k$-module which is isomorphic to the irreducible $P_k$-module $V^{(\lambda^{(k)})}$ indexed by $\lambda^{(k)} \in \hat{G}_k$. It follows that, for any choice of the isomorphism $\Phi\colon\mathbb{C}G \to P_n$ in (1.6), we must have that

$$z_{k,j}v_L = \Phi(z_{k,j})v_L = c_{k,j}(\lambda^{(k)})v_L. \qquad \blacksquare$$

*Remark 1.10.* The set $Z = \bigcup_{k=0}^{n} Z_k$ is a set of elements of the group algebra $\mathbb{C}G$ that all commute with each other. They generate a commutative subalgebra $T$ of $\mathbb{C}G$. The subalgebra $T$ acts diagonally on the basis $v_L$, i.e. for each $t \in T$ and each $L \in \mathcal{L}$,

$$tv_L = c(t,L)v_L$$

for some constant $c(t,L) \in \mathbb{C}$.

For each $\mu \in \hat{G}_k$ let $c_k(\mu)$ be the ordered $r_k$-tuple $c_k(\mu) = (c(z_{k,j}, \mu))_{1 \le j \le r_k}$. Define the *weight* of a path $L = (\lambda^{(0)}, \dots, \lambda^{(n)})$ in $\Gamma$ to be the $n$-tuple

$$\mathrm{wt}(L) = \big(c_0(\lambda^{(0)}), \cdots, c_n(\lambda^{(n)})\big). \tag{1.11}$$

The following proposition shows that in many cases the isomorphism $\Phi$ in (1.6) can be determined more or less explicitly.

**Proposition 1.12.** *The choice of an isomorphism $\Phi$ as in (1.6) is determined by the choice of elements*

$$\Phi(E_{ML}) = e_{ML} \in \mathbb{C}G$$

*for each pair of paths $L, M$ in $\Gamma$. Assume that each path in $\Gamma$ is distinguished by its weight, i.e., if $L$ and $M$ are paths in $\Gamma$ and $L \ne M$, then $\mathrm{wt}(L) \ne \mathrm{wt}(M)$.*

    (a) *For each path $L$ in $\Gamma$ the element $e_{LL}$ is determined uniquely by the elements $z_{k,j} \in Z_k$ and the constants $c_{k,j}(\mu)$, $\mu \in \hat{G}_k$, $0 \le k \le n$.*

    (b) *If $M$ and $L$ are paths in $\Gamma$ such that $M \ne L$ then $e_{ML}$ is determined up to a constant by the elements $z_{k,j} \in Z_k$ and the constants $c_{k,j}(\mu)$, $\mu \in \hat{G}_k$, $0 \le k \le n$.*

*Proof.* Let $L = (\lambda^{(0}, \dots, \lambda^{(n)})$ be a path in $\Gamma$. For each $0 \le k \le n$ and each $1 \le j \le r_k$, let

$$p_{k,j}(\lambda^{(k)}) = \prod_{c_{k,j}(\mu) \ne c_{k,j}(\lambda^{(k)})} \frac{z_{k,j} - c_{k,j}(\mu)}{c_{k,j}(\lambda^{(k)}) - c_{k,j}(\mu)}$$

where the product is over all $c_{k,j}(\mu)$, $\mu \in \hat{G}_k$ such that $c_{k,j}(\mu) \ne c_{k,j}(\lambda^{(k)})$. There may be elements $\mu \in \hat{G}_k$ such that $\mu \ne \lambda^{(k)}$ but such that $c_{k,j}(\mu) = c_{k,j}(\lambda^{(k)})$. These $\mu \in \hat{G}_k$ are not included in the product. It follows from Lemma (1.9) that if $M = (\mu^{(0)} \to \cdots \to \mu^{(n)})$ is a path in $\Gamma$ then, for any isomorphism $\Phi$ as in (1.6),

$$\Phi(p_{k,j}(\lambda^{(k)}))v_M = \begin{cases} v_M, & \text{if } c_{k,j}(\mu^{(k)}) = c_{k,j}(\lambda^{(k)}), \\ 0, & \text{otherwise.} \end{cases}$$

Define

$$e_{LL} = \prod_{k,j} p_{k,j}(\lambda^{(k)}).$$

If $M = (\mu^{(0)} \to \cdots \to \mu^{(n)})$ is a path in $\Gamma$ then

$$\Phi(e_{LL})v_M = \delta_{LM}v_M = E_{LL}v_M,$$

since, if $L \neq M$ then $\mathrm{wt}(L) \neq \mathrm{wt}(M)$. The result follows since $\Phi$ is injective.

(b) Assume that $M$ and $L$ are paths in $\Gamma$ such that $M \neq L$. Let $a \in \mathbb{C}G$ such that $e_{MM}ae_{LL} \neq 0$. Then $e_{ML}$ must be a constant times the element $e_{MM}ae_{LL} \in \mathbb{C}G$. Since the elements $e_{LL}$ and $e_{MM}$ are completely determined by the elements $z_{k,j}$ and the constants $c_{k,j}(\mu)$ it follows that the elements $e_{ML}$ are determined (up to a constant) by them. ∎

*Remark 1.13.*   Suppose that

$$\Phi: \quad \begin{matrix} P_n & \to & \mathbb{C}G \\ E_{ML} & \mapsto & e_{ML} \end{matrix} \qquad \text{and} \qquad \Phi': \quad \begin{matrix} P_n & \to & \mathbb{C}G \\ E_{ML} & \mapsto & e'_{ML} \end{matrix}$$

are two isomorphisms between the path algebra $P_n$ and $\mathbb{C}G$. Let $\kappa_{ML} \in \mathbb{C}$ be such that $e_{ML} = \kappa_{ML}e'_{ML}$. The constants $\kappa_{ML}$ must satisfy the relations

$$\kappa_{ML}\kappa_{LM} = 1, \qquad \text{and} \qquad \kappa_{ML}\kappa_{LN} = \kappa_{MN},$$

for all choices of paths $M, L, N \in \mathcal{L}$. These relations follow from the relations $e_{ML}e_{PS} = \delta_{LP}e_{MS}$ in (1.2).

*Example 1.14.*   Suppose that $\{1\} = G_0 \subseteq \cdots \subseteq G_n = G$ is a chain of finite groups such that, for each $1 \leq i \leq n$, the restriction rules describing the decomposition of irreducible $G_i$-representations into irreducible $G_{i-1}$-representations are multiplicity free. For each $1 \leq i \leq n$, let $Z_i$ be the set of sums of elements in each conjugacy class of $G_i$. Clearly $Z_i$ is a set of central elements in $\mathbb{C}G_i$. This is an example of a situation in which the paths in the graph $\Gamma$ are distinguished by their weights.

## 2. Weyl groups and Iwahori-Hecke algebras

The branching rules for the chains of Weyl groups

$$S_1 \subseteq S_2 \subseteq \cdots \subseteq S_n,$$
$$WB_2 \subseteq WB_3 \subseteq \cdots \subseteq WB_n,$$
$$WB_2 \subseteq WB_3 \subseteq WF_4,$$
$$WD_5 \subseteq WE_6 \subseteq WE_7,$$

are all multiplicity free. Thus, the Weyl groups $S_n$, $WB_n$, $WF_4$, $WE_6$, $WE_7$, all fall into the situation of example (1.14) and one can use the sets $Z_k$ consisting of all conjugacy class sums and Proposition (1.12) to compute all the irreducible representations of these Weyl groups and their corresponding Iwahori-Hecke algebras (In the $WE_7 \subseteq WE_8$ case the branching rule has multiplicities $\leq 2$ and could be treated in a similar fashion to Example 1.14 except that one would have to also use some additional elements from the centralizer of the $WE_7$ action on irreducible $WE_8$ representations.) We shall show in the remainder of this paper that one can use much smaller sets for the $Z_k$ and obtain the same results in a quicker way. We shall obtain "seminormal" representations of the Weyl groups and Iwahori-Hecke algebras of type $D_n$ by using the representation theory for type $B_n$, see section 5.

*Weyl groups*

Let $\emptyset = R_0 \subseteq \cdots \subseteq R_n = R$ be a chain of root systems and let

$$\{1\} = W_0 \subseteq \cdots \subseteq W_n$$

be the corresponding chain of Weyl groups. Let $R_k^+$ denote the set of positive roots in the root system $R_k$ and, for each $\alpha \in R_k^+$ let $s_\alpha$ denote the element of $W_k$ which is the reflection in the hyperplane perpendicular to $\alpha$.

(2.1a) If all roots in $R_k$ are the same length then the set of elements $\{s_\alpha \mid \alpha \in R_k^+\}$, is a conjugacy class in $W_k$. It follows that

$$z_{k,\ell} = \sum_{\alpha \in R_k^+} s_\alpha$$

is a central element of $\mathbb{C}W_k$. (The index $\ell$ here simply denotes that this is a sum over the "long" roots in $R_k^+$.)

(2.1b) If the roots in $R_k$ are not all the same length then there are two lengths of roots. Let $R_{k,s}^+$ be the set of short positive roots and let $R_{k,\ell}^+$ be the set of long positive roots in $R_k$. The sets $\{s_\alpha \mid \alpha \in R_{k,s}^+\}$ and $\{s_\alpha \mid \alpha \in R_{k,\ell}^+\}$ are conjugacy classes in $W$ and the elements

$$z_{k,s} = \sum_{\alpha \in R_{k,s}^+} s_\alpha, \qquad \text{and} \qquad z_{k,\ell} = \sum_{\alpha \in R_{k,\ell}^+} s_\alpha,$$

are central elements in $\mathbb{C}W_k$.

(2.1c) If the longest element $w_{k,0}$ in the Weyl group $W_k$ acts as $-1$ in the reflection representation of $W_k$ then the element

$$z_{k,0} = w_{k,0}$$

is central in $W_k$.

For each $0 \le k \le n$, let $Z_k$ denote the set of central elements in $\mathbb{C}W_k$ which are determined by (2.1). Depending on which cases apply, the set $Z_k$ contains either 1, 2, or 3 elements. In view of Lemma (1.9), we define, for each irreducible character $\chi$ of the Weyl group $W_k$ and each $z_{k,j} \in Z_k$, a constant

$$c_{k,j}(\chi) = \frac{\chi(z_{k,j})}{\chi(1)}. \tag{2.2}$$

We shall use the central elements in the sets $Z_k$ (with some slight modification in the $D_n$ case) to compute seminormal representations for the Weyl groups of types $A_{n-1}$, $B_n$, $D_n$, $G_2$.

*Remark 2.3.* In my view, the central elements in (2.1) are the appropriate generalization of Jucys-Murphy elements to arbitrary Weyl groups (or Coxeter groups). Remarks (3.5) and (4.6) illustrate this idea in special cases.

*Iwahori-Hecke algebras*

Let $\Delta$ be a Dynkin diagram of a finite Weyl group, let $R$ be the corresponding root system, and let $W$ be the corresponding Weyl group. Let $\{\alpha_i \mid i \in \Delta\}$ be the set of simple roots in $R$, indexed by the nodes in the Dynkin diagram $\Delta$. The Iwahori-Hecke algebra $H(p^2, q^2)$ corresponding to the Weyl group $W$ is the algebra over $\mathbb{C}(p, q)$ generated by elements $T_i$, $i \in \Delta$, and relations

(a) $\underbrace{T_i T_j T_i T_j \cdots}_{m_{ij} \text{ factors}} = \underbrace{T_j T_i T_j T_i \cdots}_{m_{ij} \text{ factors}}$,

where $m_{ij}$ is the order of the element $s_{\alpha_i} s_{\alpha_j}$ in the Weyl group $W$,

(b) $T_i^2 = \begin{cases} (p - p^{-1})T_i + 1, & \text{if } \alpha_i \text{ is a short root,} \\ (q - q^{-1})T_i + 1, & \text{if } \alpha_i \text{ is a long root.} \end{cases}$

If all roots in $R$ are the same length then we make the convention, for the purposes of the definition of the Iwahori-Hecke algebra, that all roots in $R$ are long and we simply define the Iwahori-Hecke algebra as an algebra $H(q^2)$ over $\mathbb{C}(q)$.

It is a standard fact ([Bou] Chap. IV §2 Ex. 23-24 or [CR]) that the Iwahori-Hecke algebra $H(q^2)$ corresponding to a Weyl group $W$ is split-semisimple and its irreducible representations can be indexed by the same set $\hat{W}$ that indexes the irreducible representations of $W$. Following the standard notation, if $w \in W$ then $T_w = T_{i_1} \cdots T_{i_p}$ where $w = s_{i_1} \cdots s_{i_p}$ is a reduced expression for $w \in W$. The element $T_w \in H$ is well defined and does not depend on the choice of the reduced expression for $w$.

Let $\emptyset = \Delta_0 \subseteq \Delta_1 \subseteq \cdots \subseteq \Delta_n = \Delta$ be a chain of Dynkin diagrams of finite Weyl groups. These Dynkin diagrams correspond to a chain of root systems $R_0 \subseteq R_1 \subseteq \cdots \subseteq R_n = R$ and to a chain of Weyl groups $W_0 \subseteq W_1 \subseteq \cdots \subseteq W_n = W$ such that, for each $1 \le i \le n$, the group $W_{i-1}$ is a parabolic subgroup of $W_i$. Let

$$\mathbb{C} \cong H_0 \subseteq \cdots \subseteq H_n = H$$

be the corresponding chain of Iwahori-Hecke algebras.

**Proposition 2.4.** *Let $H_k$ be the Iwahori-Hecke algebra corresponding to a finite Weyl group $W_k$. Let $w_{k,0}$ denote the longest element in the Weyl group $W_k$.*

*(a) The element $T_{w_{k,0}}^2$ is central in $H_k(p^2, q^2)$. If $\rho$ is an irreducible representation of $H_k(p^2, q^2)$ corresponding to the irreducible character $\chi$ of the Weyl group $W_k$, then*

$$\rho(T_{w_{k,0}}^2) = p^{2c_{k,s}(\chi)} q^{2c_{k,\ell}(\chi)} \mathrm{Id},$$

*where $c_{k,s}(\chi)$ and $c_{k,\ell}(\chi)$ are the constants given in (2.2).*

*(b) If $w_{k,0} = -1$ in the reflection representation of $W_k$, then $T_{w_{k,0}}$ is a central element in $H_k(p^2, q^2)$. If $\rho$ is an irreducible representation of $H_k(p^2, q^2)$ corresponding to the irreducible character $\chi$ of the Weyl group $W_k$, then*

$$\rho(T_{w_{k,0}}) = c_{k,0}(\chi) p^{c_{k,s}(\chi)} q^{c_{k,\ell}(\chi)} \mathrm{Id},$$

*where $c_{k,s}(\chi)$, $c_{k,\ell}(\chi)$, and $c_{k,0}(\chi)$ are the constants given in (2.2).*

*Proof.* (a) By a theorem of Breiskorn-Saito [BS] and Deligne [De], the element $T_{w_{k,0}}^2$ is central in the generalized braid group. Thus $T_{w_{k,0}}^2$ is central in $H_k(p^2, q^2)$ and it follows that $T_{w_{k,0}}^2$ acts by a constant in every irreducible representation. The constant is computed by writing $T_{w_{k,0}}$ as a product of generators and taking the determinant of both sides of the equation $T_{w_{k,0}}^2 = (T_{i_1} \cdots T_{i_N})^2 = p^{c_1} q^{c_2} \mathrm{Id}$. It remains only to note that the number of short (resp. long) roots in $R_k$ is the same as the number of factors $T_{i_j}$ in the product $T_{w_{k,0}} = T_{i_1} \cdots T_{i_N}$ such that $\alpha_{i_j}$ is a short (resp. long) root in $R_k$.

(b) The result of Brieskorn-Saito and Deligne says that $T_{w_{k,0}}$ is central in the braid group when $w_{k,0} = -1$ in the reflection representation of the Weyl group $W_k$. It follows that $T_{w_{k,0}}$ is central in $H_k(p^2, q^2)$. The eigenvalues of $T_{w_{k,0}}$ must be square roots of the eigenvalues of $T_{w_{k,0}}^2$ and they must specialize to the eigenvalues of $w_{k,0}$ when $p = q = 1$. The result now follows from (a), (2.1)(c), and the definition of the constant $c_{k,0}(\chi)$. ∎

*Remark 2.5.* The above proposition is somewhat folklore in the subject of Iwahori-Hecke algebras. The argument given here appears in Propositions 26 and 27 of Kilmoyer's thesis [Ki] and also appears in complete detail in the recent paper of Geck and Michel [GM].

Let $w_{k,0}$ be the longest element in the Weyl group $W_k$ and define $Z_k = \{z_k\}$ where

$$z_k = T_{w_{k,0}}, \qquad \text{if } w_{k,0} = -1 \text{ in the reflection representation of } W_i, \text{ and}$$
$$z_k = T_{w_{k,0}}^2, \qquad \text{otherwise.}$$

$$(2.6)$$

We shall use these central elements (with some slight modification in the $D_n$ case) to compute seminormal representations for the Iwahori-Hecke algebras of types $A_{n-1}$, $B_n$, $D_n$, $G_2$. The cases $F_4$, $E_6$ and $E_7$ will be treated in a future work.

*Remark 2.7.* In my view, the central elements in Proposition (2.4) are the analogues of the Jucys-Murphy elements for the Iwahori-Hecke algebras. Remarks (3.16) and (4.22) illustrate this idea in special cases.

The following proposition describes concretely the connection between the central elements of $H_k(p^2, q^2)$ in Proposition (2.4) and the central elements in $\mathbb{C}W_k$ given in (2.1).

**Proposition 2.8.** *If $x \in H_k(p^2, q^2)$ use the notation $[x]_{q=1}$ to denote the*

*value of $x$ when $q$ is specialized to 1. Then*

$$\left[\frac{([T_{w_{k,0}}]_{p=1})^2 - 1}{q - q^{-1}}\right]_{q=1} = \sum_{\alpha \in R_{k,s}^+} s_\alpha = z_{k,s},$$

$$\left[\frac{([T_{w_{k,0}}]_{q=1})^2 - 1}{p - p^{-1}}\right]_{p=1} = \sum_{\alpha \in R_{k,\ell}^+} s_\alpha = z_{k,\ell},$$

$$[T_{w_{k,0}}]_{p=q=1} = w_{k,0} = z_{k,0},$$

*where $z_{k,s}$, $z_{k,\ell}$, and $z_{k,0}$ are the central elements of $\mathbb{C}W_k$ given in (2.1).*

*Proof.* Let us assume, for convenience, that all roots in $R_k$ are the same length. In this case we have only one indeterminate $q$ and we are working in the Iwahori-Hecke algebra $H_k(q^2)$. The proof is similar in the general case.

If $T_{i_1} \cdots T_{i_N}$ is a reduced expression for $T_{w_{k,0}}$ then so is $T_{i_N} \cdots T_{i_1}$. We can expand $(T_{w_{k,0}})^2$ by using the relation $T_i^2 = (q - q^{-1})T_i + 1$ to obtain

$$\begin{aligned}
T_{w_{k,0}}^2 &= T_{i_1} \cdots T_{i_N} T_{i_N} \cdots T_{i_1} \\
&= 1 + (q - q^{-1}) \sum_{j=1}^N T_{i_1} \cdots T_{i_{j-1}} T_{i_j} T_{i_{j-1}} \cdots T_{i_1} \\
&\quad + \text{ terms divisible by } (q - q^{-1})^2.
\end{aligned}$$

It follows that

$$\left[\frac{T_{w_{k,0}}^2 - 1}{q - q^{-1}}\right]_{q=1} = \sum_{j=1}^N s_{i_1} \cdots s_{i_{j-1}} s_{i_j} s_{i_{j-1}} \cdots s_{i_1}.$$

The result now follows from [Bou] Chapt. VI §1, Cor 2. ∎

*Remark 2.9.* In Proposition (2.8) we have been carefree about the process of specializing $p$ and $q$ to 1. Of course this really should be done properly. One must define a $\mathbb{Z}$-form of the Iwahori-Hecke algebra $H_k(p^2, q^2)$ as an algebra over $\mathcal{A} = \mathbb{Z}[q, q^{-1}, p, p^{-1}]$ and only specialize, by an appropriate tensor product $\mathbb{Z} \otimes_{\mathcal{A}} H_k(p^2, q^2)$, elements $x$ which are in the $\mathbb{Z}$-form of $H_k(p^2, q^2)$. This is standard and it is clear that the elements $(T_{w_{k,0}}^2 - 1)/(q - q^{-1})$ in (the proof of) Proposition (2.8) are elements in the $\mathbb{Z}$-form of $H_k(q^2)$.

## 3. Type $A_{n-1}$, the symmetric group $S_n$

*The Weyl group*

The Weyl group of the root system $A_{n-1}$ is the symmetric group $S_n$ of permutations of $\{1, 2 \ldots, n\}$. The simple transpositions

$$s_i = (i-1, i), \qquad 2 \leq i \leq n,$$

generate $S_n$ and these elements satisfy the relations

$$
\begin{aligned}
s_i s_j &= s_j s_i, & |i-j| &> 1, \\
s_i s_{i+1} s_i &= s_{i+1} s_i s_{i+1}, & 2 \leq i &\leq n-1, \\
s_i^2 &= 1, & 2 \leq i &\leq n.
\end{aligned}
\tag{3.1}
$$

*Partitions and Standard Tableaux*

As in [Mac], we shall identify each partition $\lambda$ with its Ferrers diagram and say that a box $b$ in $\lambda$ is in position $(i, j)$ in $\lambda$ if $b$ is in row $i$ and column $j$ of $\lambda$. The rows and columns of $\lambda$ are labeled in the same way as for matrices. We shall write $|\lambda| = n$ if $\lambda$ is a partition with $n$ boxes. We shall often refer to partitions as *shapes*.

A *standard tableau* $L$ of shape $\lambda$ is a filling of the Ferrers diagram of $\lambda$ with the numbers $1, 2, \ldots, n$ such that the numbers are increasing left to right across the rows of $L$ and increasing down the columns of $L$. For any shape $\lambda$, let $\mathcal{L}(\lambda)$ denote the set of standard tableaux of shape $\lambda$ and, for each standard tableau $L$, let $L(k)$ denote the box containing $k$ in $L$. For example,

| 1 | 2 | 4 |
|---|---|---|
| 3 | 5 | 6 |
| 7 | 8 |   |

Figure (3.2)

is a standard tableau of shape (332).

*The chain $A_0 \subseteq A_1 \subseteq \cdots \subseteq A_{n-1}$*

of the Weyl group $S_k$ does not act by $-1$ in the reflection representation. For each $1 \le k \le n$, the set $Z_k$ contains a single element $z_{k,\ell}$, which is the central element of $\mathbb{C}S_k$ given by

$$z_{k,\ell} = \sum_{\alpha \in A_{k-1}^+} s_\alpha = \sum_{1 \le i < j \le k} (i,j) \quad . \tag{3.4}$$

(The $\ell$ here is spurious, as in (2.1a) it only indicates the fact that in the root system $A_{k-1}$ all roots are long). Since the elements $z_{k,\ell}$, $1 \le k \le n$ all commute with each other in $\mathbb{C}S_n$, it follows that the elements

$$m_k = z_{k,\ell} - z_{k-1,\ell} = \sum_{i=2}^{k} (i-1,k), \qquad 2 \le k \le n, \tag{3.5}$$

all commute with each other in $\mathbb{C}S_n$.

*Remark 3.6.* The elements $m_k$, $2 \le k \le n$, are the elements defined by Jucys [Ju1-3] and Murphy [Mu1-3].

*Weights*

The *content* of a box $b$ in a shape $\lambda$ is given by

$$\mathrm{ct}(b) = j - i, \qquad \text{if } b \text{ is in position } (i,j) \text{ in } \lambda. \tag{3.7}$$

It follows immediately from [Mac] I §7 Ex. 7 and [Mac] I §1 Ex.3 that, for each $1 \le k \le n$ and for each partition $\mu$ such that $|\mu| = k$,

$$c_{k,\ell}(\mu) = \chi^\mu(z_{k,\ell})/\chi^\mu(1) = \sum_{b \in \mu} \mathrm{ct}(b), \tag{3.8}$$

where $\chi^\mu$ denotes the character of the irreducible representation of the symmetric group $S_k$ labeled by the partition $\mu$. Following (1.11), the *weight* of a standard tableau $L = (\lambda^{(1)} \to \cdots \to \lambda^{(n)})$, where $|\lambda^{(k)}| = k$, is

$$\mathrm{wt}(L) = \big(c_{1,\ell}(\lambda^{(1)}), \ldots, c_{n,\ell}(\lambda^{(n)})\big).$$

Note that $\mathrm{wt}(L)$ is completely determined by the $n$-tuple

$$\widetilde{\mathrm{wt}}(L) = \big(\mathrm{ct}(L(1)), \ldots, \mathrm{ct}(L(n))\big), \qquad \text{since} \quad c_{k,\ell}(\lambda^{(k)}) = \sum_{i=1}^{k} \mathrm{ct}(L(i)).$$

**Proposition 3.9.** *Each standard tableau* $L = (\lambda^{(1)} \to \cdots \to \lambda^{(n)})$ *is determined uniquely by its weight.*

*Proof.* Two boxes $b$ and $b'$ in a partition $\lambda$ have the same content only if they lie on the same diagonal. It follows easily that, if $\lambda^{(i)}$ is a partition, then each of the boxes $b$ that can be added to $\lambda^{(i)}$ to get a new partition has a different content $\mathrm{ct}(b)$. Thus, the shape $\lambda^{(i+1)}$ in a standard tableau $L$ is completely determined by the previous shape $\lambda^{(i)}$ and the content $\mathrm{ct}(b)$ of the added box $b$. It follows that a standard tableau $L$ is completely determined $\widetilde{\mathrm{wt}}(L)$ and therefore by its weight $\mathrm{wt}(L)$. ∎

Proposition (1.12) and Proposition (3.9) together show that the seminormal representations of $S_n$ corresponding to the chain of groups $\{1\} = S_0 \subseteq \cdots \subseteq S_n$ are essentially determined by the elements $z_{k,\ell}$ in (3.4) and the constants $c_{k,\ell}(\mu)$ in (3.8). It follows that we should be able to determine seminormal representations of the group $S_n$ from the elements $m_k$ and the constants $\mathrm{ct}(b)$. This is done in Theorems (3.12) and (3.14) below.

*Seminormal representations*

Let $P_1 \subseteq P_2 \subseteq \cdots \subseteq P_n$ be the path algebras, defined in (1.2), which are associated to the diagram $\Gamma$ which describes the restriction rules for the chain $S_1 \subseteq \cdots \subseteq S_n$. For each partition $\lambda$ of size $n$, let

$$V^\lambda = \mathbb{C}\text{-}span\{v_L \mid L \in \mathcal{L}(\lambda)\}, \qquad (3.10)$$

so that the vectors $v_L$, indexed by standard tableaux $L$ of shape $\lambda$, form a seminormal basis of the $P_n$-module $V^\lambda$. It follows from Lemma (1.9), that for any choice of an isomorphism $\Phi$ between the path algebra $P_n$ and $\mathbb{C}S_n$ such that $\Phi(P_k) = \mathbb{C}S_k \subseteq \mathbb{C}S_n$ for all $1 \le k \le n$, we have that

$$z_{k,\ell} v_L = c_{k,\ell}(\lambda^{(k)}) v_L,$$

if $L = (\lambda^{(1)} \to \cdots \to \lambda^{(n)})$. If $m_k$ is as in (3.5), then

$$m_k v_L = \mathrm{ct}(L(k)) v_L,$$

for a standard tableau $L = (\lambda^{(1)} \to \cdots \to \lambda^{(n)})$.

For each $2 \le k \le n$ and each standard tableau $L$ of size $n$, define

$$(s_k)_{LL} = \frac{1}{\mathrm{ct}(L(k)) - \mathrm{ct}(L(k-1))} \qquad (3.11)$$

In the interests of space we shall not give the proof of the following theorems here. The proofs are essentially the same as the proofs which are given for Theorem (4.15) and Theorem (4.18).

**Theorem 3.12.** (Young [Y]) *Let $\lambda$ be a partition such that $|\lambda| = n$. Define an action of each generator $s_2, \ldots, s_n$ of the symmetric group $S_n$ on $V^\lambda$ by defining*

$$s_i v_L = (s_i)_{LL} v_L \; + \; (1 + (s_i)_{LL}) v_{s_i L}, \qquad 2 \leq i \leq n, \qquad (3.13)$$

*where $s_i L$ is the same standard tableau as $L$ except that the positions of $i$ and $i-1$ are switched in $s_i L$. If $s_i L$ is not standard, then we define $v_{s_i L} = 0$. This action extends to a well defined action of $S_n$ on $V^\lambda$.*

**Theorem 3.14.** (Young [Y]) *The $S_n$ modules $V^\lambda$ defined in Theorem (3.12), where $\lambda$ runs over all partitions such that $|\lambda| = n$, form a complete set of nonisomorphic irreducible modules for the symmetric group $S_n$ and, for each $\lambda$, the basis $\{v_L \mid L \in \mathcal{L}(\lambda)\}$ is a seminormal basis of the $S_n$-module $V^\lambda$.*

## Iwahori-Hecke algebras $HA_{n-1}(q^2)$

Let $q$ be an indeterminate. The Iwahori-Hecke algebra $HA_{k-1}(q^2)$ corresponding to the root system $A_{k-1}$ is the associative algebra with 1 over the field $\mathbb{C}(q)$ given by generators $T_2, T_3, \ldots, T_k$ and relations

$$\begin{aligned}
T_i T_j &= T_j T_i, \qquad |i - j| > 1, \\
T_i T_{i+1} T_i &= T_{i+1} T_i T_{i+1}, \qquad 2 \leq i \leq k-1, \\
T_i^2 &= (q - q^{-1}) T_i + 1, \qquad 2 \leq i \leq k.
\end{aligned} \qquad (3.15)$$

*Analogues of Jucys-Murphy elements*

For each $2 \leq k \leq n$, define

$$M_k = T_k \cdots T_3 T_2 T_2 T_3 \cdots T_k. \qquad (3.16)$$

In type $A_{k-1}$ the longest element $w_{k,0}$ of the group $S_k$ does not act by $-1$ in the reflection representation. Following (2.4), we define sets $Z_k = \{z_k\}$, $2 \leq k \leq n$, where $z_k$ is the central element of $HA_{k-1}(q^2)$ given by

$$z_k = T_{w_{k,0}}^2 = M_k M_{k-1} \cdots M_2.$$

Since the elements $z_k$, $2 \le k \le n$, all commute in $HA_{n-1}(q^2)$ it follows that the elements $M_k$, $2 \le k \le n$ all commute with each other.

*Remark 3.17.*  Using the relation $T_i^2 = (q - q^{-1})T_i + 1$, an easy computation shows that

$$\frac{M_k - 1}{q - q^{-1}} = \sum_{i=2}^{k} T_{(i-1, k)}, \qquad (3.18)$$

where $T_{(i-1, k)} = T_k T_{k-1} \cdots T_{i+1} T_i T_{i+1} \cdots T_k$.  The elements in (3.18) are elements used by Dipper, James, and Murphy in their work on Iwahori-Hecke algebras of type $A$, see [DJ1-2], [M1], [M4].  It is clear that (3.18) gives a $q$-analogue of the Jucys-Murphy elements in (3.5).

*Seminormal representations*

Let $P_1 \subseteq P_2 \subseteq \cdots \subseteq P_n$ be the path algebras (over the field $\mathbb{C}(q)$ instead of $\mathbb{C}$), defined in (1.2), which are associated to the diagram $\Gamma$ which describes the restriction rules for the chain $S_1 \subseteq \cdots \subseteq S_n$.  For each partition $\lambda$ of size $n$, let

$$V^\lambda = \mathbb{C}(q)\text{-}span\{v_L \mid L \in \mathcal{L}(\lambda)\}, \qquad (3.19)$$

so that the vectors $v_L$, indexed by standard tableaux $L$ of shape $\lambda$, form a seminormal basis of the $P_n$-module $V^\lambda$.  It follows from Lemma (1.9), that for any choice of an isomorphism $\Phi$ between the path algebra $P_n$ and $HA_{n-1}(q^2)$ such that $\Phi(P_k) = HA_{k-1}(q^2) \subseteq HA_{n-1}(q^2)$ for all $1 \le k \le n$, we have that

$$z_k v_L = T_{w_{k,0}}^2 v_L = q^{c_{k,\ell}(\lambda^{(k)})} v_L,$$

if $L = (\lambda^{(1)} \subseteq \cdots \subseteq \lambda^{(n)})$ and $c_{k,\ell}(\lambda^{(k)})$ is as given in (3.8).  Thus,

$$M_k v_L = T_k \cdots T_3 T_2 T_2 T_3 \cdots T_k v_L = T_{w_{k,0}}^2 T_{w_{k-1,0}}^{-2} = q^{2\operatorname{ct}(L(k))} v_L, \qquad (3.20)$$

if $L = (\lambda^{(1)} \to \cdots \to \lambda^{(n)})$ is a standard tableau.  For each $2 \le k \le n$ and each standard tableau $L$ of size $n$, define

$$(T_k)_{LL} = \frac{q - q^{-1}}{1 - \frac{\operatorname{CT}(L(k-1))}{\operatorname{CT}(L(k))}} \quad \text{where} \quad \operatorname{CT}(b) = q^{2\operatorname{ct}(L(k))}, \qquad (3.21)$$

In the interests of space we shall not give the proof of the following theorems here.  The proofs are essentially the same as the proofs which are given for Theorem (4.26) and Theorem (4.28).

**Theorem 3.22.** (Hoefsmit [H]) *Let $\lambda$ be a partition such that $|\lambda| = n$. Define an action of each generator $T_2, \ldots, T_n$ of the Iwahori-Hecke algebra $HA_{n-1}(q^2)$ on $V^\lambda$ by defining*

$$T_i v_L = (T_i)_{LL} v_L \ + \ (q^{-1} + (T_i)_{LL}) v_{s_i L}, \qquad 2 \le i \le n,$$

*where $s_i L$ is the same standard tableau as $L$ except that the positions of $i$ and $i-1$ are switched in $s_i L$. If $s_i L$ is not standard, then we define $v_{s_i L} = 0$. This action extends to a well defined action of $HA_{n-1}(q^2)$ on $V^\lambda$.*

**Theorem 3.23.** (Hoefsmit [H]) *The $HA_{n-1}(q^2)$ modules $V^\lambda$ defined in Theorem (3.22), where $\lambda$ runs over all partitions such that $|\lambda| = n$, form a complete set of nonisomorphic irreducible modules for the Iwahori-Hecke algebra $HA_{n-1}(q^2)$ and, for each $\lambda$, the basis $\{v_L \mid L \in \mathcal{L}(\lambda)\}$ is a seminormal basis of the $HA_{n-1}(q^2)$-module $V^\lambda$.*

## 4. Type $B_n$, $n \geq 2$

*The Weyl group*

The Weyl group $WB_n$ of type $B_n$ is the group of signed permutations of $1, 2, \ldots, n$. More specifically, $WB_n$ consists of all permutations $\pi$ of $\{-n, \ldots, -1, 1, \ldots, n\}$ such that $\pi(-k) = -\pi(k)$ for all $1 \leq k \leq n$. We represent elements of $WB_n$ in cycle notation as permutations of $\{-n, \ldots, -1, 1, \ldots, n\}$. The elements

$$s_1 = (1, -1), \qquad \text{and} \qquad s_i = (i-1, i)(-(i-1), -i), \quad 2 \leq i \leq n,$$

generate $WB_n$ and satisfy the relations

$$\begin{aligned}
s_i s_j &= s_j s_i, \qquad |i-j| > 1, \\
s_i s_{i+1} s_i &= s_{i+1} s_i s_{i+1}, \qquad 2 \leq i \leq n-1, \\
s_1 s_2 s_1 s_2 &= s_2 s_1 s_2 s_1, \\
s_i^2 &= 1, \qquad 2 \leq i \leq n.
\end{aligned} \tag{4.1}$$

*Double Partitions and Standard Tableaux*

A *double partition* of size $n$, $\lambda = (\alpha, \beta)$, is an ordered pair of partitions $\alpha$ and $\beta$ such that $|\alpha| + |\beta| = n$. We shall often refer to double partitions as *shapes*. A *standard tableau* $L = (L^\alpha, L^\beta)$ of shape $\lambda = (\alpha, \beta)$ is a filling of the Ferrers diagram of $\lambda$ with the numbers $1, 2, \ldots, n$ such that the numbers are increasing left to right across the rows of $L^\alpha$ and $L^\beta$ and increasing down the columns of $L^\alpha$ and $L^\beta$. For any shape $\lambda$, let $\mathcal{L}(\lambda)$ denote the set of standard tableaux of shape $\lambda$ and, for each standard tableau $L$, let $L(k)$ denote the box containing $k$ in $L$. For example,

Figure (4.2)

is a standard tableau of shape $((332),(411))$.

*The chain $B_0 \subseteq B_1 \subseteq \cdots \subseteq B_n$*

By convention we let $B_0 = \emptyset$ be the empty root system and $B_1 = A_1$. The chain of root systems $B_0 \subseteq B_1 \subseteq \cdots \subseteq B_n$ corresponds to the chain of Weyl groups

$$\{1\} \subseteq WB_1 \subseteq WB_2 \subseteq \cdots \subseteq WB_n, \qquad (4.3)$$

where $WB_k$ denotes the hyperoctahedral group of signed permutations of $1, 2, \ldots, k$.

The irreducible representations of the symmetric group $WB_k$ are indexed by double partitions $\lambda = (\alpha, \beta)$ such that $|\lambda| = |\alpha| + |\beta| = k$. The restriction rule from $WB_k$ to $WB_{k-1}$ is given by

$$V^{(\alpha,\beta)} \Big\downarrow_{WB_{k-1}}^{WB_k} \cong \bigoplus_{(\mu,\nu)\in(\alpha,\beta)^-} V^{(\mu,\nu)}$$

where the sum is over all double partitions $(\mu, \nu)$ of size $k-1$ that are gotten from $(\alpha, \beta)$ by removing one box. If we define the graph $\Gamma$ as in (1.1) for the chain in (4.3), then a path $(\lambda^{(0)} \to \cdots \to \lambda^{(n)})$ in $\Gamma$ is naturally identified with the standard tableau $L$ of shape $\lambda^{(n)}$ which has $i$ in the box which is added to obtain $\lambda^{(i)}$ from $\lambda^{(i-1)}$. The graph $\Gamma$ for the case of the chain $\{1\} \subseteq WB_1 \subseteq WB_2 \subseteq WB_3$, is displayed in Figure (4.4).

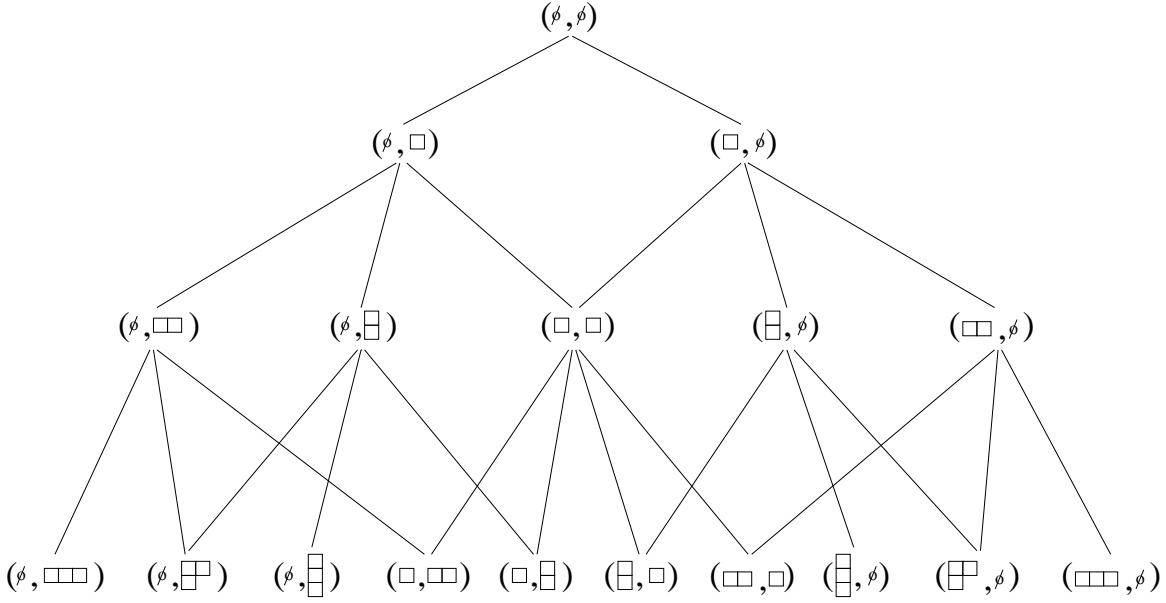

Figure (4.4)

*Analogues of Jucys-Murphy elements*

Following (2.1), let us compute the sets $Z_k$ for this case. In the root system $B_k$, $k \geq 2$, we have both long and short roots and the longest element $w_{k,0} = (1,-1)(2,-2)\cdots(k,-k)$ of the Weyl group $WB_k$ acts by $-1$ in the reflection representation. For each $1 \leq k \leq n$, let

$$Z_k = \{z_{k,s}, z_{k,\ell}, z_{k,0}\},$$

where

$$
\begin{aligned}
z_{k,s} &= \sum_{\alpha \in (B_k)_s^+} s_\alpha = \sum_{i=1}^{k} (i,-i), \\
z_{k,\ell} &= \sum_{\alpha \in (B_k)_\ell^+} s_\alpha = \sum_{1 \leq i < j \leq k} (i,j) + (i,-j)(-i,j), \\
z_{k,0} &= w_{k,0} = (1,-1)(2,-2)\cdots(k,-k).
\end{aligned}
\tag{4.5}
$$

In (4.5) the set $(B_k)_s^+$ (resp. $(B_k)_\ell^+$) is the set of short (resp. long) positive roots in the root system $B_k$. Since the elements $z_{k,j}$, $j \in \{s, \ell, 0\}$, $1 \leq k \leq n$ all commute with each other in $\mathbb{C}WB_n$, it follows that the elements

$$m_{k,s} = z_{k,s} - z_{k-1,s} = (k, -k), \quad 1 \le k \le n,$$

$$m_{k,\ell} = z_{k,\ell} - z_{k-1,\ell} = \sum_{i=2}^{k} (i-1, k) + (i-1, -k)(-(i-1), k), \quad 2 \le k \le n,$$

all commute with each other in $\mathbb{C}WB_n$.

*Remark 4.6.* The elements $m_{k,s}$ and $m_{k,\ell}$ are the appropriate $B_n$-analogues of the Jucys-Murphy elements (3.5) for the symmetric group. Cherednik [Ch] has used a linear combination of $m_{k,s}$ and $m_{k,\ell}$ as an analogue of the Jucys-Murphy element.

*Weights*

The *sign* and the *content* of a box $b$ in a shape $(\alpha, \beta)$ are given respectively by

$$\operatorname{sgn}(b) = \begin{cases} 1, & \text{if } b \in \alpha, \\ -1, & \text{if } b \in \beta, \end{cases} \tag{4.7}$$

$$\operatorname{ct}(b) = j - i, \qquad \text{if } b \text{ is in position } (i, j) \text{ (in either } \alpha \text{ or } \beta).$$

**Proposition 4.8.** *Fix $1 \le k \le n$ and let $z_{k,s}$, $z_{k,\ell}$, and $z_{k,0}$ be the central elements in $\mathbb{C}WB_k$ given in (4.5). Let $(\mu, \nu)$ be a double partition such that $|\mu| + |\nu| = k$. Let $\chi^{(\mu,\nu)}$ be the character of the irreducible representation of $WB_k$ indexed by the double partition $(\mu, \nu)$. Then, in the notation of Lemma (1.9) and (4.7),*

$$c_{k,s}(\mu, \nu) = \chi^{(\mu,\nu)}(z_{k,s})/\chi^{(\mu,\nu)}(1) = \sum_{b \in (\mu,\nu)} \operatorname{sgn}(b),$$

$$c_{k,\ell}(\mu, \nu) = \chi^{(\mu,\nu)}(z_{k,\ell})/\chi^{(\mu,\nu)}(1) = \sum_{b \in (\mu,\nu)} \operatorname{ct}(b),$$

$$c_{k,0}(\mu, \nu) = \chi^{(\mu,\nu)}(z_{k,0})/\chi^{(\mu,\nu)}(1) = \prod_{b \in (\mu,\nu)} \operatorname{sgn}(b).$$

*Proof.* Fix $k$ and a double partition $(\mu, \nu)$ such that $|\mu| + |\nu| = k$. Let $a = |\mu|$ and $b = |\nu|$. Let $WB_a$ be the subgroup of $WB_k$ of signed permutations of $\{1, 2, \ldots, a\}$ and let $WB_b$ be the subgroup of $WB_k$ of signed permutations of $\{a+1, \ldots, k\}$. Let $S_a$ and $S_b$ be the symmetric groups of permutations of $\{1, \ldots, a\}$ and $\{a+1, \ldots, k\}$ respectively.

Let $V^\mu$ be the irreducible module for the symmetric group $S_a$ which is labeled by the partition $\mu$. Extend this module to the group $WB_a$ by letting the signs, i.e. the element $(1, -1)$ act trivially. Let $V^\nu$ be the irreducible module for the symmetric group $S_b$ which is labeled by $\nu$. Extend this module to group $WB_b$ by letting the element $(a+1, -(a+1))$ act by $-1$ on $V^\nu$. Then $V^\mu$ and $V^\nu$ are irreducible modules for the groups $WB_a$ and $WB_b$ respectively, they are the modules that are ordinarily denoted by $V^{(\mu, \emptyset)}$ and $V^{(\emptyset, \nu)}$ respectively. It follows from this construction of $V^{(\mu, \emptyset)}$ and $V^{(\emptyset, \nu)}$ that

$$
\begin{aligned}
&\chi^{(\mu, \emptyset)}((1, 2)) = \chi^\mu((1, 2)), \\
&\chi^{(\mu, \emptyset)}((1, -1)) = \chi^\mu(1), \\
&\chi^{(\mu, \emptyset)}((1, -1)(2, -2) \cdots (a, -a)) = \chi^\mu(1), \\
&\chi^{(\emptyset, \nu)}((a+1, a+2)) = \chi^\nu((a+1, a+2)), \\
&\chi^{(\emptyset, \nu)}((a+1, -(a+1)) = -\chi^\nu(1)), \\
&\chi^{(\emptyset, \nu)}((a+1, -(a+1)(a+2, (a+2)) \cdots (k, -k)) = (-1)^b \chi^\nu(1)),
\end{aligned}
\tag{4.9}
$$

where $\chi^\mu$ and $\chi^\nu$ denote the irreducible characters of the symmetric groups $S_a$ and $S_b$ labeled by $\mu$ and $\nu$, respectively.

It is well known that the induced module

$$
V^{(\mu, \emptyset)} \otimes V^{(\emptyset, \nu)} \uparrow_{WB_a \times WB_b}^{WB_k} \cong V^{(\mu, \nu)}
$$

is a realization of the irreducible $WB_k$-module indexed by the double partition $(\mu, \nu)$. (I believe that this fact is originally due to Specht [Spe].) Given this realization we can write down its character $\chi^{(\mu, \nu)}$ explicitly by using the standard formula for induced characters.

$$
\chi^{(\mu, \nu)}(w) = \sum_{g_i^{-1} w g_i \in WB_a \times WB_b} \chi^{(\mu, \emptyset)}(g_i^{-1} w g_i) \chi^{(\emptyset, \nu)}(g_i^{-1} w g_i),
\tag{4.10}
$$

where the sum is over coset representatives $g_i$ of $WB_k/(WB_a \times WB_b)$ such that $g_i^{-1} w g_i \in WB_a \times WB_b$.

Using (4.9) and (4.10),

$$\frac{\chi^{(\mu,\nu)}\left(\sum_{\alpha\in(B_k)_\ell^+} s_\alpha\right)}{\chi^{(\mu,\nu)}(1)}$$

$$= \frac{k(k-1)\chi^{(\mu,\nu)}((1,2))}{\chi^{(\mu,\nu)}(1)}$$

$$= \frac{k(k-1)\left(\binom{k-2}{a-2}\chi^{(\mu,\emptyset)}((1,2))\chi^{(\emptyset,\nu)}(1) + \binom{k-2}{a}\chi^{(\mu,\emptyset)}(1)\chi^{(\emptyset,\nu)}((a+1,a+2))\right)}{\binom{k}{a}\chi^{(\mu,\emptyset)}(1)\chi^{(\emptyset,\nu)}(1)}$$

$$= k(k-1)\left(\frac{a(a-1)}{k(k-1)}\frac{\chi^\mu((1,2))}{\chi^\mu(1)} + \frac{b(b-1)}{k(k-1)}\frac{\chi^\nu((a+1,a+2))}{\chi^\nu(1)}\right)$$

$$= 2\left(\frac{\chi^\mu\left(\sum_{\alpha\in A_{a-1}^+} s_\alpha\right)}{\chi^\mu(1)} + \frac{\chi^\nu\left(\sum_{\alpha\in A_{b-1}^+} s_\alpha\right)}{\chi^\nu(1)}\right).$$

The value for $c_{k,\ell}(\mu,\nu)$ is obtained from this last equation and (3.8). A similar calculation gives

$$\frac{\chi^{(\mu,\nu)}\left(\sum_{\alpha\in(B_k)_s^+} s_\alpha\right)}{\chi^{(\mu,\nu)}(1)}$$

$$= \frac{k\chi^{(\mu,\nu)}((1,-1))}{\chi^{(\mu,\nu)}(1)}$$

$$= \frac{k\left(\binom{k-1}{a-1}\chi^{(\mu,\emptyset)}((1,-1))\chi^{(\emptyset,\nu)}(1) + \binom{k-1}{a}\chi^{(\mu,\emptyset)}(1)\chi^{(\emptyset,\nu)}((a+1,-(a+1)))\right)}{\binom{k}{a}\chi^{(\mu,\emptyset)}(1)\chi^{(\emptyset,\nu)}(1)}$$

$$= \frac{a\chi^\mu((1,-1))}{\chi^\mu(1)} + \frac{b\chi^\nu((a+1,-(a+1)))}{\chi^\nu(1)}$$

$$= a - b.$$

Since $\text{sgn}(b) = 1$ for all boxes $b \in \mu$, $\text{sgn}(b) = -1$ for all boxes $b \in \nu$, and $|\mu| = a$ and $|\nu| = b$ the formula for $c_{k,s}(\mu,\nu)$ follows. The formula for $c_{k,0}(\mu,\nu)$ is obtained in a similar fashion. ∎

Following (1.11), the *weight* of a standard tableau $L = (\lambda^{(1)} \to \cdots \to \lambda^{(n)})$, where $|\lambda^{(k)}| = k$, is

$$\text{wt}(L) = \left(c_1(\lambda^{(1)}),\ldots,c_n(\lambda^{(n)})\right),$$

where, for a double partition $(\mu, \nu)$, $c_k(\mu, \nu)$ is the triple

$$c_k(\mu, \nu) = (c_{k,s}(\mu, \nu), c_{k,\ell}(\mu, \nu), c_{k,0}(\mu, \nu)),$$

determined by Proposition (4.8). Note that $\mathrm{wt}(L)$ is completely determined by the $n$-tuples

$$\widetilde{\mathrm{wt}}_1(L) = \big(\mathrm{ct}(L(1)), \ldots, \mathrm{ct}(L(n))\big), \text{ and } \widetilde{\mathrm{wt}}_2(L) = \big(\mathrm{sgn}(L(1)), \ldots, \mathrm{sgn}(L(n))\big),$$

since $c_{k,s}(\lambda^{(k)}) = \sum_{i=1}^{k} \mathrm{sgn}(L(i))$, $c_{k,\ell}(\lambda^{(k)}) = \sum_{i=1}^{k} \mathrm{ct}(L(i))$, and $c_{k,0}(\lambda^{(k)}) = \prod_{i=1}^{k} \mathrm{sgn}(L(i))$.

**Proposition 4.11.** *Each standard tableau $L = (\lambda^{(1)} \to \cdots \to \lambda^{(n)})$ is determined uniquely by its weight.*

*Proof.* Suppose that the weight $\mathrm{wt}(L)$ of a standard tableau $L$ is given but that we don't know $L$. The vector $\mathrm{wt}(L)$ uniquely determines the vectors $\widetilde{\mathrm{wt}}_1(L)$ and $\widetilde{\mathrm{wt}}_2(L)$. We want to show that the tableau $L$ can be reconstructed from $\widetilde{\mathrm{wt}}_1(L)$ and $\widetilde{\mathrm{wt}}_2(L)$. Assume that we have reconstructed $L$ up to the $i$th step, i.e., assume that we know $\lambda^{(0)}, \ldots, \lambda^{(i)}$, but that we don't yet know $\lambda^{(i+1)}$. Suppose that $\lambda^{(i)}$ is the double partition $(\mu^{(i)}, \nu^{(i)})$.

We need to figure out from $\widetilde{\mathrm{wt}}_1(L)$ and $\widetilde{\mathrm{wt}}_2(L)$ where to add the box to get $\lambda^{(i+1)} = (\mu^{(i+1)}, \nu^{(i+1)})$. The entry $\mathrm{sgn}(L(i+1))$ in $\widetilde{\mathrm{wt}}_2(L)$ tells us whether we must add the box to the partition $\mu^{(i)}$ or to the partition $\nu^{(i)}$. As in the proof of Proposition (3.9), the entry $\mathrm{ct}(L(i+1))$ from $\widetilde{\mathrm{wt}}_1(L)$ indicates the position where this box must be added. It follows that $\lambda^{(i+1)}$ is uniquely determined. Thus $L$ is completely determined by $\mathrm{wt}(L)$. ∎

Proposition (1.12) and Proposition (4.11) together show that the seminormal representations of $WB_n$ corresponding to the chain of groups $\{1\} = WB_0 \subseteq \cdots \subseteq WB_n$ are essentially determined by the elements $z_{k,j}$ in (4.5) and the constants $c_{k,j}(\mu)$ in Proposition (4.8). It follows that we should be able to determine seminormal representations of the group $WB_n$ from the elements $m_{k,j}$ and the constants $\mathrm{ct}(b)$ and $\mathrm{sgn}(b)$. This is done in Theorems (4.15) and (4.18) below.

*Seminormal representations*

Let $P_1 \subseteq P_2 \subseteq \cdots \subseteq P_n$ be the path algebras, defined in (1.2), which are associated to the diagram $\Gamma$ which describes the restriction rules for the chain $WB_1 \subseteq \cdots \subseteq WB_n$. For each double partition $(\alpha, \beta)$ such that $|\alpha| + |\beta| = n$, let

$$V^{(\alpha,\beta)} = \mathbb{C}\text{-}span\{v_L \mid L \in \mathcal{L}(\alpha, \beta)\}, \tag{4.12}$$

so that the vectors $v_L$, indexed by standard tableaux $L$ of shape $(\alpha, \beta)$, form a seminormal basis of the $P_n$-module $V^{(\alpha,\beta)}$. It follows from Lemma (1.9), that for any choice of an isomorphism $\Phi$ between the path algebra $P_n$ and $\mathbb{C}WB_n$ such that $\Phi(P_k) = \mathbb{C}WB_k \subseteq \mathbb{C}WB_n$ for all $1 \le k \le n$,

$$z_{k,s} v_L = c_{k,s}(\lambda^{(k)}) v_L, \quad z_{k,\ell} v_L = c_{k,\ell}(\lambda^{(k)}) v_L, \quad z_{k,0} v_L = c_{k,0}(\lambda^{(k)}) v_L,$$

if $L = (\lambda^{(1)} \to \cdots \to \lambda^{(n)})$. Thus, by Proposition (4.8),

$$m_{k,s} v_L = \text{sgn}(L(k)) v_L, \quad \text{and} \quad m_{k,\ell} v_L = \text{ct}(L(k)) v_L, \tag{4.13}$$

if $L = (\lambda^{(1)} \to \cdots \to \lambda^{(n)})$.

For each $2 \le k \le n$ and each standard tableau $L$ of length $n$, define

$$(s_k)_{LL} = \frac{1 + \text{sgn}(L(k))\text{sgn}(L(k-1))}{\text{ct}(L(k)) - \text{ct}(L(k-1))} \tag{4.14}$$

**Theorem 4.15.** (Young [Y]) *Let $(\alpha, \beta)$ be a double partition such that $|\alpha| + |\beta| = n$. Define an action of each generator $s_1, \ldots, s_n$ of $WB_n$ on $V^{(\alpha,\beta)}$ by defining*

$$\begin{aligned} s_1 v_L &= \text{sgn}(L(1)) v_L, \\ s_i v_L &= (s_i)_{LL} v_L \ + \ (1 + (s_i)_{LL}) v_{s_i L}, \qquad 2 \le i \le n, \end{aligned} \tag{4.16}$$

*where $s_i L$ is the same standard tableau as $L$ except that the positions of $i$ and $i-1$ are switched in $s_i L$. If $s_i L$ is not standard, then we define $v_{s_i L} = 0$. This action extends to a well defined action of $WB_n$ on $V^{(\alpha,\beta)}$.*

*Proof.* We shall show that the action of $s_i$, $1 \le i \le n$, on $V^{(\alpha,\beta)}$ is essentially forced by the formulas in (4.13). We shall prove this for $i = n$. The proof for $i < n$ is similar. Note that the formula for the action of $s_1$ follows immediately from the formula for the action of $m_{1,s}$ in (4.13).

For any two standard tableaux $M$ and $L$ of shape $(\alpha, \beta)$ let $(s_n)_{ML}$ be the coefficient of $v_M$ in $s_n v_L$.

*Step 1.* Let $L = (\lambda^{(0)} \to \cdots \to \lambda^{(n)})$ be a standard tableau of shape $(\alpha, \beta)$. For each $0 \le k \le n$ and each $j \in \{s, \ell, 0\}$ $1 \le j \le r_k$, let

$$p_{k,j}(\lambda^{(k)}) = \prod_{c_{k,j}(\mu) \neq c_{k,j}(\lambda^{(k)})} \frac{z_{k,j} - c_{k,j}(\mu)}{c_{k,j}(\lambda^{(k)}) - c_{k,j}(\mu)},$$

as in the proof of Proposition (1.12). Define

$$p_{L[n-2]} = \prod_{k=1}^{n-2} \prod_{j \in \{s, \ell, 0\}} p_{k,j}(\lambda^{(k)}). \tag{$*$}$$

If $M = (\mu^{(0)} \to \cdots \to \mu^{(n)})$ is another standard tableau of shape $(\alpha, \beta)$ then

$$p_{L[n-2]} v_M = \begin{cases} v_M, & \text{if } \mu^{(k)} = \lambda^{(k)}, \text{ for } 1 \le k \le n-2, \\ 0, & \text{otherwise.} \end{cases}$$

Note that since $\lambda^{(n-2)}$ and $\lambda^{(n)} = (\alpha, \beta)$ only differ by two boxes there are only two tableaux $M$ such that $\mu^{(k)} = \lambda^{(k)}$ for all $1 \le k \le n-2$. These two tableaux are $s_n L$ and $L$ itself. It follows that

$$p_{L[n-2]} v_M = \begin{cases} v_M, & \text{if } M = s_n L \text{ or } M = L, \\ 0, & \text{otherwise.} \end{cases}$$

Since each of the elements $z_{k,s}$, $z_{k,\ell}$, $z_{k,0}$ appearing the product (*) is an element of $WB_{n-2}$ it follows that $p_{L[n-2]}$ commutes with $s_n$ in $WB_n$. Thus

$$(s_n)_{LL} v_L + (s_n)_{s_n L, L} v_{s_n L} = p_{L[n-2]} s_n v_L = s_n p_{L[n-2]} v_L = s_n v_L.$$

It follows that $(s_n)_{ML} = 0$ unless $M = s_n L$ or $M = L$.

*Step 2.* A direct computation shows that

$$s_n m_{n-1,\ell} = m_{n,\ell} s_n - 1 - m_{n,s} m_{n-1,s}.$$

Let both sides act on $v_L$ and take the coefficient of $v_L$ in the result. Then, using (4.13), we have

$$(s_n)_{LL} \text{ct}(L(n-1)) = \text{ct}(L(n))(s_n)_{LL} - 1 - \text{sgn}(L(n))\text{sgn}(L(n-1)).$$

It follows that $(s_n)_{LL}$ is as given in (4.14).

*Step 3.* Consider the equation $s_n^2 = 1$. Let both sides act on $v_L$ and take the coefficient of $v_L$ in the result. We get the equation $(s_n)_{LL}^2 + (s_n)_{LM}(s_n)_{ML} = 1$, where $M = s_n L$. It follows that

$$(s_n)_{LM}(s_n)_{ML} = (1 + (s_n)_{LL})(1 - (s_n)_{LL}). \tag{4.17}$$

By Proposition (1.12)(b), the values of $(s_n)_{LM}$ and $(s_n)_{ML}$ are determined only up to a constant and we may choose them to be anything such that the equation in (4.17) holds. Note that this is consistent with the definition of the action in the statement of the theorem since $1 - (s_n)_{s_n L, s_n L} = 1 + (s_n)_{LL}$. ∎

**Theorem 4.18.** (Young [Y]) *The $WB_n$ modules $V^{(\alpha,\beta)}$ defined in Theorem (4.15), where $(\alpha, \beta)$ runs over all ordered pairs of partitions such that $|\alpha| + |\beta| = n$, form a complete set of nonisomorphic irreducible modules for the Weyl group $WB_n$ and, for each $(\alpha, \beta)$, the basis $\{v_L \mid L \in \mathcal{L}(\alpha, \beta)\}$ is a seminormal basis of the $WB_n$-module $V^{(\alpha,\beta)}$.*

*Proof.* This now follows immediately by induction on $n$. Indeed, $V^{(\alpha,\beta)}$ is the unique $WB_n$-module such that
(1) The equations for the action of $z_{n,s}$, $z_{n,\ell}$, $z_{n,0}$ are as in (4.13), and
(2) On restriction to $WB_{n-1}$ we have that

$$V^{(\alpha,\beta)} \big\downarrow_{WB_{n-1}}^{WB_n} \cong \bigoplus_{(\mu,\nu) \in (\alpha,\beta)^-} V^{(\mu,\nu)}$$

where the sum is over all double partitions $(\mu, \nu)$ of size $n - 1$ that are gotten from $(\alpha, \beta)$ by removing one box. ∎

## Iwahori-Hecke algebras $HB_n(p^2, q^2)$

Let $p$ and $q$ be indeterminates. The Iwahori-Hecke algebra $HB_k(p^2, q^2)$ corresponding to the root system $WB_k$ is the associative algebra with 1 over the field $\mathbb{C}(p, q)$ given by generators $T_1, T_2, \ldots, T_k$ and relations

$$
\begin{aligned}
T_{s_i} T_{s_j} &= T_{s_j} T_{s_i}, &&|i - j| > 1, \\
T_{s_i} T_{s_{i+1}} T_{s_i} &= T_{s_{i+1}} T_{s_i} T_{s_{i+1}}, &&2 \leq i \leq n-1, \\
T_{s_1} T_{s_2} T_{s_1} T_{s_2} &= T_{s_2} T_{s_1} T_{s_2} T_{s_1}, && \\
T_{s_1}^2 &= (p - p^{-1}) T_1 + 1, && \\
T_{s_i}^2 &= (q - q^{-1}) T_{s_i} + 1, &&1 \leq i \leq n.
\end{aligned}
\tag{4.19}
$$

*Analogues of Jucys-Murphy elements*

For each $1 \leq k \leq n$, define

$$M_k = T_k \cdots T_2 T_1 T_2 \cdots T_k. \qquad (4.20)$$

The longest element $w_{k,0}$ in the Weyl group $WB_k$ acts by $-1$ in the reflection representation. Following (1.8) and Proposition (2.4), we define sets $Z_k = \{z_k\}$, $1 \leq k \leq n$, where $z_k$ is the central element of $HB_k(p^2, q^2)$ given by

$$z_k = T_{w_{k,0}} = M_k M_{k-1} \cdots M_2 M_1. \qquad (4.21)$$

Since the elements $z_k$, $1 \leq k \leq n$, all commute in $HB_n(p^2, q^2)$ it follows that the elements $M_k$, $1 \leq k \leq n$ all commute with each other.

*Remark 4.22.* The elements $M_k$ appear in Hoefsmit [H] Proposition 3.3.3 and also in work of Ariki-Koike [AK], Ariki [Ar], Broué-Malle [BM], and Dipper-James [DJ3]. These elements can be viewed as the quantized versions of the elements in (4.13).

*Seminormal representations*

Let $P_1 \subseteq P_2 \subseteq \cdots \subseteq P_n$ be the path algebras (over $\mathbb{C}(p,q)$ instead of $\mathbb{C}$), defined in (1.2), which are associated to the diagram $\Gamma$ which describes the restriction rules for the chain $WB_1 \subseteq \cdots \subseteq WB_n$. For each double partition $(\alpha, \beta)$ such that $|\alpha| + |\beta| = n$, let

$$V^{(\alpha, \beta)} = \mathbb{C}(p,q) \text{-}span\{v_L \mid L \in \mathcal{L}(\alpha, \beta)\}, \qquad (4.23)$$

so that the vectors $v_L$, indexed by standard tableaux $L$ of shape $(\alpha, \beta)$, form a seminormal basis of the $P_n$-module $V^{(\alpha, \beta)}$. It follows from Lemma (1.9), that for any choice of an isomorphism $\Phi$ between the path algebra $P_n$ and $HB_n(p^2, q^2)$ such that $\Phi(P_k) = HB_k(p^2, q^2) \subseteq HB_n(p^2, q^2)$ for all $1 \leq k \leq n$,

$$z_k v_L = T_{w_{k,0}} v_L = c_{k,0}(\lambda^{(k)}) p^{c_{k,s}(\lambda^{(k)})} q^{c_{k,\ell}(\lambda^{(k)})} v_L,$$

if $L = (\lambda^{(1)} \to \cdots \to \lambda^{(n)})$. Thus, by Proposition (4.8),

$$\begin{aligned}
M_k v_L &= T_{k-1} \cdots T_2 T_1 T_2 \cdots T_{k-1} v_L = T_{w_{k,0}} T_{w_{k-1,0}}^{-1} v_L \\
&= \text{sgn}(L(k)) p^{\text{sgn}(L(k))} q^{2\text{ct}(L(k))} v_L,
\end{aligned} \qquad (4.24)$$

if $L = (\lambda^{(1)} \to \cdots \to \lambda^{(n)})$ is a standard tableau. For each $2 \leq k \leq n$ and each standard tableau $L$ of size $n$, define

$$(T_k)_{LL} = \frac{q - q^{-1}}{1 - \frac{\text{CT}(L(k-1))}{\text{CT}(L(k))}} \quad \text{where} \quad \text{CT}(b) = \text{sgn}(L(k)) p^{\text{sgn}(L(k))} q^{2\text{ct}(L(k))}, \tag{4.25}$$

if $b$ is a box in a shape $\lambda = (\alpha, \beta)$.

**Theorem 4.26.** *Let $(\alpha, \beta)$ be a double partition such that $|\alpha| + |\beta| = n$. Define an action of each generator $T_1, \ldots, T_n$ of $HB_n(p^2, q^2)$ on $V^{(\alpha,\beta)}$ by defining*

$$\begin{aligned}
T_1 v_L &= \text{CT}(L(1)) v_L, \\
T_i v_L &= (T_i)_{LL} v_L \ + \ (q^{-1} + (T_i)_{LL}) v_{s_i L}, \qquad 2 \leq i \leq n,
\end{aligned} \tag{4.27}$$

*where $s_i L$ is the same standard tableau as $L$ except that the positions of $i$ and $i-1$ are switched in $s_i L$. If $s_i L$ is not standard, then we define $v_{s_i L} = 0$. This action extends to a well defined action of $HB_n(p^2, q^2)$ on $V^{(\alpha,\beta)}$.*

*Proof.* The proof is similar to the proof of Theorem (4.15), in all essential aspects. We shall only give the details for step 2.

*Step 2.* It is immediate from the definition of $M_k$ in (4.20) that $M_n = T_n M_{n-1} T_n$, which can be rewritten as

$$T_n^{-1} = M_n^{-1} T_n M_{n-1}.$$

Rewrite $T_n^{-1}$ as $T_n - (q - q^{-1})$, let both sides act on $v_L$ and take the coefficient of $v_L$ in the result. Using (4.24), we get

$$(T_n)_{LL} - (q - q^{-1}) = \text{CT}(L(n))^{-1} (T_n)_{LL} \text{CT}(L(n-1)).$$

It follows that $(T_n)_{LL}$ is as given in (4.25). ■

As in the case of Weyl group $WB_n$, Theorem (4.18), the following result follows easily.

**Theorem 4.28.** (Hoefsmit [H], Thm. 2.2.14) *The $HB_n(p^2, q^2)$-modules $V^{(\alpha,\beta)}$ defined in Theorem (4.26), where $(\alpha, \beta)$ runs over all ordered pairs of partitions such that $|\alpha| + |\beta| = n$, form a complete set of nonisomorphic irreducible modules for the Iwahori-Hecke algebra $HB_n(p^2, q^2)$. For each $(\alpha, \beta)$, the basis $\{v_L \mid L \in \mathcal{L}(\alpha, \beta)\}$ is a seminormal basis of the $HB_n(p^2, q^2)$-module $V^{(\alpha,\beta)}$.*

## 5. Type $D_n$, $n \geq 4$

*The Weyl group*

The Weyl group $WD_n$ of type $D_n$ is the group of signed permutations of $\{1, 2, \ldots, n\}$ with an even number of signs. More specifically, $WD_n$ consists of all permutations $\pi$ of $\{-n, \ldots, -1, 1, \ldots, n\}$ such that $\pi(-k) = -\pi(k)$ for all $1 \leq k \leq n$, and an even number of the elements of $\{\pi(1), \pi(2), \ldots, \pi(n)\}$ are negative. We represent elements of $WD_n$ in cycle notation as permutations of $\{-n, \ldots, -1, 1, \ldots, n\}$.

The elements

$$\tilde{s}_1 = (1, -2)(2, -1), \qquad \text{and} \qquad \tilde{s}_i = (i-1, i)(-(i-1), -i), \quad 2 \leq i \leq n,$$

generate $WD_n$ and satisfy the relations

$$\begin{aligned}
\tilde{s}_i \tilde{s}_j &= \tilde{s}_j \tilde{s}_i, & |i - j| > 1, \text{ and } i, j > 1, \\
\tilde{s}_1 \tilde{s}_j &= \tilde{s}_j \tilde{s}_1, & \text{if } j \neq 3, \\
\tilde{s}_1 \tilde{s}_3 \tilde{s}_1 &= \tilde{s}_3 \tilde{s}_1 \tilde{s}_3, & \\
\tilde{s}_i \tilde{s}_{i+1} \tilde{s}_i &= \tilde{s}_{i+1} \tilde{s}_i \tilde{s}_{i+1}, & 2 \leq i \leq n-1, \\
\tilde{s}_i^2 &= 1, & 1 \leq i \leq n.
\end{aligned} \tag{5.1}$$

The Weyl group $WD_n$ can be realized as a normal subgroup of the Weyl group $WB_n$ of index 2 by defining

$$\tilde{s}_1 = s_1 s_2 s_1, \qquad \text{and} \qquad \tilde{s}_i = s_i, \quad 2 \leq i \leq n, \tag{5.2}$$

where $s_i$, $1 \leq i \leq n$, are as in (4.1).

*Double Partitions and Standard Tableaux*

We shall use the same notations for partitions, double partitions, shapes, and tableaux as in section 4. For each standard tableau $L = (L^\alpha, L^\beta)$ of shape $(\alpha, \beta)$ define $\sigma L$ to be the standard tableau of shape $(\beta, \alpha)$ given by $\sigma L = (L^\beta, L^\alpha)$,

$$\begin{aligned}
\sigma: \quad \mathcal{L}(\alpha, \beta) &\rightarrow \mathcal{L}(\beta, \alpha) \\
(L^\alpha, L^\beta) &\mapsto (L^\beta, L^\alpha).
\end{aligned} \tag{5.3}$$

The map $\sigma$ is an involution on the set of standard tableaux whose shape is a double partition.

*Which chain?*

One finds that it is more natural to use the representation theory of the Weyl groups $WB_n$ and the fact that $WD_n$ is a normal subgroup of index 2 in $WB_n$ rather than to try to choose an appropriate chain of root systems leading up to $D_n$. The reason for this is that one wants to have an approach that treats all of the groups $WD_n$, $n \geq 4$, uniformly. Otherwise one must distinguish the cases when $n$ is even and when $n$ is odd. In the end we shall find a set of commuting elements in the group algebra of $WD_n$, analogues of the Jucys-Murphy elements, which determine a complete set of irreducible representations.

*Representations*

We shall retain the notations from section 4 for the sign and the content of a box in a double partition. Let $\lambda = (\alpha, \beta)$ be a double partition such that $|\alpha| + |\beta| = n$. As in (4.12), let

$$V^{(\alpha,\beta)} = \mathbb{C}\text{-}span\{v_L \mid L \in \mathcal{L}(\alpha, \beta)\} \tag{5.4}$$

so that the vectors $v_L$ form a basis of the vector space $V^{(\alpha,\beta)}$ indexed by standard tableaux $L$ of shape $(\alpha, \beta)$.

For each standard tableau $L$, define

$$(s_k)_{LL} = \frac{1 + \mathrm{sgn}(L(k))\mathrm{sgn}(L(k-1))}{\mathrm{ct}(L(k)) - \mathrm{ct}(L(k-1))}, \qquad \text{for } 2 \leq k \leq n, \tag{5.5}$$

as in (4.14). Recall, Theorem (4.15), that there is an action of $WB_n$ on the vector space $V^{(\alpha,\beta)}$. Restricting this action to $WD_n$ gives

$$\begin{aligned}
\tilde{s}_1 v_L &= s_1 s_2 s_1 v_L = (s_2)_{LL} v_L - (1 + (s_2)_{LL}) v_{s_2 L}, \\
\tilde{s}_i v_L &= s_i v_L = (s_i)_{LL} v_L + (1 + (s_i)_{LL}) v_{s_i L}, \qquad 2 \leq i \leq n,
\end{aligned} \tag{5.6}$$

for each $L \in \mathcal{L}(\alpha, \beta)$, where, as in the case of type $B_n$, we define $v_{s_i L} = 0$ if $s_i L$ is not standard. In deriving the first formula of (5.6) it is helpful to note that $\mathrm{sgn}(L(1)) = \pm 1$, and $\mathrm{sgn}(s_2 L(1)) = -\mathrm{sgn}(L(1))$ if $s_2 L$ is standard.

Now suppose $n$ is even, and let $\alpha$ be a partition such that $2|\alpha| = n$. Define

$$\begin{aligned}
V^{(\alpha,\alpha)^+} &= \mathbb{C}\text{-}span\{v_L + v_{\sigma L} \mid L \in \mathcal{L}(\alpha, \alpha)\} \subseteq V^{(\alpha,\alpha)}, \\
V^{(\alpha,\alpha)^-} &= \mathbb{C}\text{-}span\{v_L - v_{\sigma L} \mid L \in \mathcal{L}(\alpha, \alpha)\} \subseteq V^{(\alpha,\alpha)}.
\end{aligned} \tag{5.7}$$

The following results (well known) follow easily from Clifford theory [CR] since $WD_n$ is a subgroup of index 2 in $WB_n$ and $\sigma$ commutes with the action of $WD_n$ on the vectors $v_L$, $L \in \mathcal{L}$.

**Proposition 5.8.**
(a) For each pair of partitions $(\alpha, \beta)$ such that $|\alpha| + |\beta| = n$, $V^{(\alpha,\beta)}$ and $V^{(\beta,\alpha)}$ are isomorphic $WD_n$-modules.
(b) For each partition $\alpha$ such that $2|\alpha| = n$, the subspaces $V^{(\alpha,\alpha)\pm}$ are $WD_n$-submodules of $V^{(\alpha,\alpha)}$, and

$$V^{(\alpha,\alpha)} \cong V^{(\alpha,\alpha)^+} \oplus V^{(\alpha,\alpha)^-},$$

as $WD_n$-modules.

**Theorem 5.9.** (Young [Y]) *The modules* $V^{(\alpha,\beta)}$, *where* $(\alpha, \beta)$ *runs over all unordered pairs of partitions such that* $\alpha \neq \beta$ *and* $|\alpha| + |\beta| = n$ *and, when* $n$ *is even, the modules* $V^{(\alpha,\alpha)^+}$ *and* $V^{(\alpha,\alpha)^-}$, *where* $\alpha$ *runs over all partitions such that* $2|\alpha| = n$, *form a complete set of nonisomorphic irreducible modules for* $WD_n$.

*Remark 5.10.* The involution $\sigma$ on standard tableaux in (5.3) is a realization of the module isomorphism between the $WD_n$-modules $V^{(\alpha,\beta)}$ and $V^{(\beta,\alpha)}$, which, in turn, comes from the automorphism of the Dynkin diagram of type $D_n$.

Instead of defining $V^{(\alpha,\alpha)\pm}$ as in (5.7), let us define them as the quotient spaces

$$V^{(\alpha,\alpha)^+} = \frac{V^{(\alpha,\alpha)}}{\langle v_L = v_{\sigma L}\rangle} \quad \text{and} \quad V^{(\alpha,\alpha)^-} = \frac{V^{(\alpha,\alpha)}}{\langle v_L = -v_{\sigma L}\rangle}, \qquad (5.11)$$

where $\sigma$ is the involution given in (5.3) and $\langle v_L = v_{\sigma L}\rangle$ and $\langle v_L = -v_{\sigma L}\rangle$ denote the subspaces spanned the the vectors $v_L - v_{\sigma L}$ and by $v_L + v_{\sigma L}$ respectively. Clearly the two definitions of the modules $V^{(\alpha,\alpha)\pm}$ are equivalent, the first represents $V^{(\alpha,\alpha)\pm}$ as subspaces of $V^{(\alpha,\alpha)}$, and the second as quotient spaces of $V^{(\alpha,\alpha)}$. The only difference is that for some computations the quotient module approach is easier, one may compute the action as in the formulas in (5.6) and then apply the relations $v_L = \pm v_{\sigma L}$.

*Analogues of Jucys-Murphy elements*

**Theorem 5.12.** *Define elements* $m_k$, $2 \leq k \leq n$, *in the group algebra of the Weyl group* $WD_n$ *by*

$$\tilde{m}_{k,1} = (1,-1)(k,-k) \quad \text{and} \quad \tilde{m}_{k,2} = \sum_{i=2}^{k}(i-1,k) + (i-1,-k)(-(i-1),k),$$

where elements of $WD_n$ are written in cycle notation as permutations of $\{-n, \ldots, -1, 1, \ldots, n\}$. Then the elements $\tilde{m}_{k,1}$ and $\tilde{m}_{k,2}$ all commute with each other in $\mathbb{C}WD_n$ and they act in the representations $V^{(\alpha, \beta)}$ and $V^{(\alpha, \alpha)\pm}$ from (5.6) and (5.7) by

$\tilde{m}_{k,1}v_L = \mathrm{sgn}(L(1))\mathrm{sgn}(L(k))v_L$, for all standard tableaux $L$, and

$\tilde{m}_{k,1}v_L^\pm = \mathrm{sgn}(L(1))\mathrm{sgn}(L(k))v_L^\pm$, for all standard tableaux $L$ of shape $(\alpha, \alpha)$,

$\tilde{m}_{k,2}v_L = \mathrm{ct}(L(k))v_L$, for all standard tableaux $L$, and

$\tilde{m}_{k,2}v_L^\pm = \mathrm{ct}(L(k))v_L^\pm$, for all standard tableaux $L$ of shape $(\alpha, \alpha)$.

*Proof.* This follows immediately from (4.13) once one notices that $\tilde{m}_{k,1} = m_{1,s}m_{k,s}$ and $\tilde{m}_{k,2} = m_{k,\ell}$, where $m_{k,s}$ and $m_{k,\ell}$ are as in (4.5). $\blacksquare$

*Weights*

If $L$ is a standard tableau of shape $(\alpha, \beta)$ define

$$\widetilde{\mathrm{wt}}_1(L) = \big(\mathrm{ct}(L(1)), \ldots, \mathrm{ct}(L(n))\big), \quad \text{and}$$
$$\widetilde{\mathrm{wt}}_2'(L) = \big(\mathrm{sgn}(L(1))^2, \ldots, \mathrm{sgn}(L(1))\mathrm{sgn}(L(n))\big),$$

where sgn and ct are as given in (4.7).

**Lemma 5.13.** *If $L$ is a standard tableau then there is only one other standard tableau $L'$ such that $\widetilde{\mathrm{wt}}_1(L') = \widetilde{\mathrm{wt}}_1(L)$ and $\widetilde{\mathrm{wt}}_2'(L') = \widetilde{\mathrm{wt}}_2'(L)$. This standard tableau is $L' = \sigma L$, where $\sigma$ is the involution defined in (5.3).*

*Proof.* Let $\widetilde{wt}_1$ and $\widetilde{wt}_2$ be as defined in Proposition (4.11). It follows from Proposition (4.11) that $\widetilde{wt}_1(L)$ and $\widetilde{wt}_2(L)$ uniquely determine $L$. Since $\widetilde{wt}_2'(L)$ is always either $+1$ or $-1$ times every entry in the sequence $\widetilde{wt}_2(L)$ it follows that there can be at most two standard tableaux $L$ and $L'$ with the same weights $\widetilde{\mathrm{wt}}_1(L') = \widetilde{\mathrm{wt}}_1(L)$ and $\widetilde{\mathrm{wt}}_2'(L') = \widetilde{\mathrm{wt}}_2'(L)$. On the other hand it is immediate that one always has that $\widetilde{\mathrm{wt}}_2'(\sigma L) = \widetilde{\mathrm{wt}}_2'(L)$ since $\mathrm{sgn}(\sigma L(k)) = -\mathrm{sgn}(L(k))$ for all $k$. $\blacksquare$

**Iwahori-Hecke algebras $HD_n(q^2)$**

Let $q$ be an indeterminate. The Iwahori-Hecke algebra $HD_n(q^2)$ of type $D_n$ is the associative algebra with 1 over the field $\mathbb{C}(q)$ given by generators $\tilde{T}_1, \tilde{T}_2, \ldots, \tilde{T}_n$ and relations

$$
\begin{aligned}
\tilde{T}_i \tilde{T}_j &= \tilde{T}_j \tilde{T}_i, && |i-j| > 1, i, j > 1, \\
\tilde{T}_1 \tilde{T}_j &= \tilde{T}_j \tilde{T}_1, && \text{if } j \neq 3, \\
\tilde{T}_1 \tilde{T}_3 \tilde{T}_1 &= \tilde{T}_3 \tilde{T}_1 \tilde{T}_3, && \\
\tilde{T}_i \tilde{T}_{i+1} \tilde{T}_i &= \tilde{T}_{i+1} \tilde{T}_i \tilde{T}_{i+1}, && 2 \leq i \leq n-1, \\
\tilde{T}_i^2 &= (q - q^{-1})\tilde{T}_i + 1, && 1 \leq i \leq n.
\end{aligned}
\tag{5.14}
$$

Let $HB_n(1, q^2)$ be the algebra over $\mathbb{C}(q)$ defined by generators $T_1, \ldots, T_n$ and relations as in (4.19) except with $p = 1$. Define

$$
\tilde{T}_1 = T_1 T_2 T_1, \qquad \text{and} \qquad \tilde{T}_i = T_i, \qquad 2 \leq i \leq n. \tag{5.15}
$$

Then one checks that with these definitions the $\tilde{T}_i$ satisfy the relations in (5.14). The elements $\tilde{T}_i$, $1 \leq i \leq n$, generate a subalgebra of the algebra $HB_n(1, q^2)$ which is isomorphic to the algebra $HD_n(q^2)$.

*Representations*

One derives the representation theory of the Iwahori-Hecke algebra $HD_n(q^2)$ using the results from section 4 and the fact that $HD_n(q^2)$ is a subalgebra of the Iwahori-Hecke algebra $HB_n(1, q^2)$. The procedure is exactly as for the case (5.4-5.9) of the Weyl groups $WD_n \subseteq WB_n$.

Let $V^{(\alpha, \beta)}$ be as in (4.23). As in (4.25), for $2 \leq k \leq n$ and each standard tableau $L$, define

$$
(\tilde{T}_k)_{LL} = (T_k)_{LL} = \frac{q - q^{-1}}{1 - \frac{\mathrm{CT}(L(k-1))}{\mathrm{CT}(L(k))}} \quad \text{where} \quad \mathrm{CT}(b) = \mathrm{sgn}(L(k)) q^{2\mathrm{ct}(L(k))}, \tag{5.16}
$$

for a box $b$ in a shape $\lambda = (\alpha, \beta)$. Restricting the action (4.27) of $HB_n(1, q^2)$ to $HD_n(q^2)$ gives

$$
\begin{aligned}
\tilde{T}_{s_1} v_L &= T_{s_1} T_{s_2} T_{s_1} v_L = (\tilde{T}_{s_2})_{LL} v_L - (q^{-1} + (\tilde{T}_{s_2})_{LL}) v_{s_2 L}, \\
\tilde{T}_{s_i} v_L &= T_{s_i} v_L = (\tilde{T}_{s_i})_{LL} v_L + (q^{-1} + (\tilde{T}_{s_i})_{LL}) v_{s_i L}, \qquad 2 \leq i \leq n,
\end{aligned}
\tag{5.17}
$$

for each $L \in \mathcal{L}(\alpha, \beta)$, where, as in the case of type $B_n$, we define $v_{s_i L} = 0$ if $s_i L$ is not standard. In deriving the first formula of (5.17) it is helpful to note that $CT(L(1)) = \pm 1$, and $CT(s_2 L(1)) = -CT(L(1))$ if $s_2 L$ is standard.

If $n$ is even, and $\alpha$ is a partition such that $2|\alpha| = n$, define $V^{(\alpha,\alpha)\pm}$ as in (5.7) except over the field $\mathbb{C}(q)$. The following results can be proved by "setting $q = 1$" and then using Proposition (5.8) and Theorem (5.9).

**Proposition 5.18.** (Hoefsmit [H], Lemmas 2.3.3 and 2.3.5)
(a) For each pair of partitions $(\alpha, \beta)$ such that $|\alpha| + |\beta| = n$, $V^{(\alpha,\beta)}$ and $V^{(\beta,\alpha)}$ are isomorphic $HD_n(q^2)$-modules.
(b) For each partition $\alpha$ such that $2|\alpha| = n$, the subspaces $V^{(\alpha,\alpha)\pm}$ are $HD_n(q^2)$-submodules of $V^{(\alpha,\alpha)}$, and

$$V^{(\alpha,\alpha)} \cong V^{(\alpha,\alpha)^+} \oplus V^{(\alpha,\alpha)^-},$$

as $HD_n(q^2)$-modules.

**Theorem 5.19.** (Hoefsmit [H], Thm. 2.3.9) The modules $V^{(\alpha,\beta)}$, where $(\alpha, \beta)$ runs over all unordered pairs of partitions such that $\alpha \neq \beta$ and $|\alpha| + |\beta| = n$ and, when $n$ is even, the modules $V^{(\alpha,\alpha)^+}$ and $V^{(\alpha,\alpha)^-}$, where $\alpha$ runs over all partitions such that $2|\alpha| = n$, form a complete set of nonisomorphic irreducible modules for $HD_n(q^2)$.

*Analogues of Jucys-Murphy elements*

Define

$$\tilde{M}_1 = 1, \quad \tilde{M}_2 = \tilde{T}_{s_2}\tilde{T}_{s_1}, \quad \text{and}$$
$$\tilde{M}_k = \tilde{T}_{s_k}\tilde{T}_{s_{k-1}} \cdots \tilde{T}_{s_3}\tilde{T}_{s_2}\tilde{T}_{s_1}\tilde{T}_{s_3}\tilde{T}_{s_4} \cdots \tilde{T}_{s_{k-1}}\tilde{T}_{s_k}, \quad \text{for } 3 \leq k \leq n. \tag{5.20}$$

If $w_0$ is the longest element of the Weyl group $WD_n$, then $\tilde{T}_{w_0} = \tilde{M}_n\tilde{M}_{n-1} \cdots \tilde{M}_1$ is the corresponding element in the Iwahori-Hecke algebra $HD_n(q^2)$.

**Theorem 5.21.** The action of the element $\tilde{M}_k$ in the irreducible representations given by Theorem (5.19) is

$$\tilde{M}_k v_L = CT(L(1))CT(L(k))v_L, \quad \text{for all standard tableaux } L, \text{ and}$$
$$\tilde{M}_k v_L^{\pm} = CT(L(1))CT(L(k))v_L^{\pm}, \quad \text{for all standard tableaux } L \text{ of shape } (\alpha, \alpha).$$

*Proof.* Let $M_k$ be the elements of $HB_n(1, q^2)$ given by (4.20), and use the imbedding of $HD_n(q^2)$ into $HB_n(1, q^2)$. The case $k = 1$ is trivial, since $CT(L(1)) = \pm 1$. For $k = 2$, observe that $\tilde{M}_2 = \tilde{T}_{s_2}\tilde{T}_{s_1} = T_{s_2}T_{s_1}T_{s_2}T_{s_1} = M_2M_1$. For $3 \le k \le n$, note that $T_{s_1}$ commutes with $T_{s_3}, T_{s_4}, \ldots$ in $HB_n(1, q^2)$, and thus

$$\begin{aligned}
\tilde{M}_k &= \tilde{T}_{s_k}\tilde{T}_{s_{k-1}} \cdots \tilde{T}_{s_3}\tilde{T}_{s_2}\tilde{T}_{s_1}\tilde{T}_{s_3}\tilde{T}_{s_4} \cdots \tilde{T}_{s_{k-1}}\tilde{T}_{s_k} \\
&= T_{s_k}T_{s_{k-1}} \cdots T_{s_3}T_{s_2}(T_{s_1}T_{s_2}T_{s_1})T_{s_3}T_{s_4} \cdots T_{s_{k-1}}T_{s_k} \\
&= T_{s_k}T_{s_{k-1}} \cdots T_{s_3}T_{s_2}T_{s_1}T_{s_2}T_{s_3}T_{s_4} \cdots T_{s_{k-1}}T_{s_k}T_{s_1} = M_kM_1.
\end{aligned}$$

The result now follows from the definition of the action of $HB_n(1, q^2)$ and of $HD_n(q^2)$ on irreducible modules and (4.24). ∎

## 6. TYPE $G_2$

*The chain $A_0 \subseteq A_1 \subseteq G_2$*

The Weyl group $WG_2$ is the dihedral group of order 12. The group $WG_2$ can be presented by generators $s_1, s_2$ and relations

$$s_1 s_2 s_1 s_2 s_1 s_2 = s_2 s_1 s_2 s_1 s_2 s_1, \quad \text{and} \quad s_i^2 = 1, \quad \text{for } i = 1, 2.$$

The irreducible representations of the dihedral group $WG_2$ can be indexed by the labels

$$\hat{G}_2 = \{\phi_{1,0}, \phi_{1,6}, \phi'_{1,3}, \phi''_{1,3}, \phi_{2,1}, \phi_{2,2}\}$$

and the character table of the group $WG_2$ is given by

|  | 1 | $s_1$ | $s_2$ | $s_1 s_2$ | $s_1 s_2 s_1 s_2$ | $s_1 s_2 s_1 s_2 s_1 s_2$ |
|---|---|---|---|---|---|---|
| $\phi_{1,0}$ | 1 | 1 | 1 | 1 | 1 | 1 |
| $\phi_{1,6}$ | 1 | $-1$ | $-1$ | 1 | 1 | 1 |
| $\phi'_{1,3}$ | 1 | 1 | $-1$ | $-1$ | 1 | $-1$ |
| $\phi''_{1,3}$ | 1 | $-1$ | 1 | $-1$ | 1 | $-1$ |
| $\phi_{2,1}$ | 2 | 0 | 0 | 1 | $-1$ | $-2$ |
| $\phi_{2,2}$ | 2 | 0 | 0 | $-1$ | $-1$ | 2 |

The chain of root systems $A_0 \subseteq A_1 \subseteq G_2$ corresponds to the chain of Weyl groups $S_1 \subseteq S_2 \subseteq WG_2$, where $S_1$ and $S_2$ are symmetric groups. The graph $\Gamma$, as defined in (1.1), corresponding to the inclusion $S_1 \subseteq S_2 \subseteq WG_2$ is given by

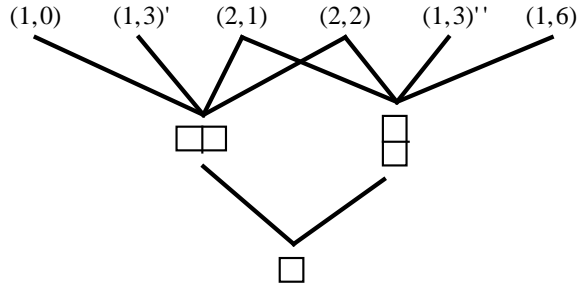

Figure (6.1)

We have indexed the representations of the symmetric groups $S_1$ and $S_2$ by partitions as in section 3.

*Analogues of Jucys-Murphy elements*

Following (1.8) and (2.1), let us compute the sets $Z_k$ for this case. In the root system $A_1$ all roots are the same length and the longest element $w_{1,0}$ in the Weyl group $S_2$ acts by $-1$ in the reflection representation. In the root system $G_2$ we have both long and short roots and the longest element $w_{2,0} = s_1 s_2 s_1 s_2 s_1 s_2$ of the Weyl group $WG_2$ acts by $-1$ in the reflection representation. Let

$$Z_1 = \{z_{1,0}\}, \quad \text{and} \quad Z_2 = \{z_{2,s}, z_{2,\ell}, z_{2,0}\},$$

where

$$\begin{aligned}
z_{1,0} &= w_{1,0} = s_1, \\
z_{2,s} &= \sum_{\alpha \in (G_2)_s^+} s_\alpha = s_1 + s_2 s_1 s_2 + s_1 s_2 s_1 s_2 s_1, \\
z_{2,\ell} &= \sum_{\alpha \in (G_2)_\ell^+} s_\alpha = s_2 + s_1 s_2 s_1 + s_2 s_1 s_2 s_1 s_2, \\
z_{2,0} &= w_{2,0} = s_1 s_2 s_1 s_2 s_1 s_2.
\end{aligned} \tag{6.2}$$

In (6.2) the set $(G_2)_s^+$ (resp. $(G_2)_\ell^+$) is the set of short (resp. long) positive roots in the root system $G_2$. The elements $z_{k,j}$ in (6.2) are the appropriate $G_2$-analogues of the Jucys-Murphy elements in (3.5).

*Weights*

Following Lemma (1.9) and (2.2), we use the character table of $WG_2$ to compute

$$c_{2,s}(\phi_{1,0}) = 3, \qquad c_{2,\ell}(\phi_{1,0}) = 3, \qquad c_{2,0}(\phi_{1,0}) = 1,$$

$$c_{2,s}(\phi_{1,6}) = -3, \qquad c_{2,\ell}(\phi_{1,6}) = -3, \qquad c_{2,0}(\phi_{1,6}) = 1,$$

$$c_{2,s}(\phi_{1,3}') = 3, \qquad c_{2,\ell}(\phi_{1,3}') = -3, \qquad c_{2,0}(\phi_{1,3}') = -1,$$

$$c_{2,s}(\phi_{1,3}'') = -3, \qquad c_{2,\ell}(\phi_{1,3}'') = 3, \qquad c_{2,0}(\phi_{1,3}'') = -1,$$

$$c_{2,s}(\phi_{2,1}) = 0, \qquad c_{2,\ell}(\phi_{2,1}) = 0, \qquad c_{2,0}(\phi_{2,1}) = -1,$$

$$c_{2,s}(\phi_{2,2}) = 0, \qquad c_{2,\ell}(\phi_{2,2}) = 0, \qquad c_{2,0}(\phi_{2,2}) = 1,$$

$$c_{1,0}((2)) = 1, \qquad c_{1,0}((1^2)) = -1,$$

so that $c_{k,j}(\mu) = \chi^\mu(z_{k,j})/\chi^\mu(1)$ where $\chi^\mu$ denotes the irreducible character labeled by $\mu$. The *weight* of a path $L = (\square \to \lambda^{(1)} \to \lambda^{(2)})$ in the graph $\Gamma$ is

$$\mathrm{wt}(L) = \left( c_{1,0}(\lambda^{(1)}), c_{2,s}(\lambda^{(2)}), c_{2,\ell}(\lambda^{(2)}), c_{2,0}(\lambda^{(2)}) \right). \tag{6.3}$$

**Proposition 6.4.** *Each path $L = (\square \to \lambda^{(1)} \to \lambda^{(2)})$ in $\Gamma$ is distinguished by its weight, i.e., if $L$ and $M$ are paths in $\Gamma$ and $L \neq M$ then $\mathrm{wt}(L) \neq \mathrm{wt}(M)$.*

*Proof.* This follows easily by a direct check. ∎

Proposition (1.12) and Proposition (6.4) together show that the seminormal representations of $WG_2$ corresponding to the chain of groups $\{1\} \subseteq S_2 \subseteq WG_2$ are essentially determined by the elements $z_{k,j}$ in (6.2) and the constants $c_{k,j}(\mu)$ which appear in (6.3). These representations are given in Theorem (6.7) below.

*Seminormal representations*

Let $P_1 \subseteq P_2$ be the path algebras, defined in section 1, which are associated to the diagram $\Gamma$ in (6.1). For each $\lambda \in \hat{G}_2$, let

$$V^\lambda = \mathbb{C}\text{-}span\{v_L \mid L \in \mathcal{L}(\lambda)\}, \tag{6.5}$$

so that the vectors $v_L$, indexed by the paths $L = (\square \to \lambda^{(1)} \to \lambda^{(2)} = \lambda)$ in $\Gamma$ which end at $\lambda$, form a seminormal basis of the irreducible $P_2$-module $V^\lambda$. It follows from Lemma (1.9), that for any choice of an isomorphism $\Phi$ between the path algebra $P_2$ and the $\mathbb{C}WG_2$ such that $\Phi(P_1) = \mathbb{C}S_2 \subseteq \mathbb{C}WG_2$,

$$\begin{aligned}
z_{1,0} v_L &= c_{1,0}(\lambda^{(1)}) v_L, & z_{2,s} v_L &= c_{2,s}(\lambda^{(2)}) v_L, \\
z_{2,0} v_L &= c_{2,0}(\lambda^{(2)}) v_L, & z_{2,\ell} v_L &= c_{2,\ell}(\lambda^{(2)}) v_L,
\end{aligned} \tag{6.6}$$

if $L = (\square \to \lambda^{(1)} \to \lambda^{(2)})$ is a path in $\Gamma$.

**Theorem 6.7.** *Irreducible seminormal representations of the Weyl group $WG_2$ with respect to the chain $S_1 \subseteq S_2 \subseteq WG_2$ can be given by*

$$\phi_{1,0}(s_1) = (1), \quad \phi_{1,6}(s_1) = (-1), \quad \phi'_{1,3}(s_1) = (1), \quad \phi''_{1,3}(s_1) = (-1),$$
$$\phi_{1,0}(s_2) = (1), \quad \phi_{1,6}(s_2) = (-1), \quad \phi'_{1,3}(s_2) = (-1), \quad \phi''_{1,3}(s_2) = (1),$$

$$\phi_{2,1}(s_1) = \begin{pmatrix} 1 & 0 \\ 0 & -1 \end{pmatrix} \qquad \phi_{2,2}(s_1) = \begin{pmatrix} 1 & 0 \\ 0 & -1 \end{pmatrix}$$

$$\phi_{2,1}(s_2) = \begin{pmatrix} \frac{1}{2} & \frac{1}{2} \\ \frac{3}{2} & -\frac{1}{2} \end{pmatrix} \qquad \phi_{2,2}(s_2) = \begin{pmatrix} -\frac{1}{2} & \frac{3}{2} \\ \frac{1}{2} & \frac{1}{2} \end{pmatrix}.$$

*Proof.* For any two paths $M$ and $L$ in the graph $\Gamma$ which end at the same label $\lambda$. Let $(s_k)_{ML}$ denote the coefficient of $v_M$ in $s_k v_L$. The constant $(s_k)_{ML}$ is a matrix entry in the matrix for $s_k$ in the irreducible representation labeled by $\lambda$. It follows from Proposition (1.12)(a) and Proposition (6.4) that the diagonal entries of these matrices must be determined by the equations in (6.6) and that the off diagonal entries are determined up to a constant.

The matrices giving the one-dimensional representations are easily gotten from the relations in (6.6). Let us explain how one derives the matrices for the two dimensional case.

(a) The matrices for $s_1$ are determined by (6.6).

(b) From the definitions (6.2), one gets easily that $z_{2,\ell} = s_2 + z_{1,\ell} s_2 z_{1,\ell} + z_{2,0} z_{1,\ell}$. Let both sides of this equation act on $v_L$ and take the coefficient of $v_L$ to get, via (6.6), the equation

$$c_{2,\ell}(\lambda^{(2)}) = (s_2)_{LL} + c_{1,\ell}(\lambda^{(1)})(s_2)_{LL} c_{1,\ell}(\lambda^{(1)}) + c_{2,0}(\lambda^{(2)}) c_{1,\ell}(\lambda^{(1)}).$$

It follows that

$$(s_2)_{LL} = \frac{c_{2,\ell}(\lambda^{(2)}) - c_{2,0}(\lambda^{(2)}) c_{1,0}(\lambda^{(1)})}{1 + c_{1,0}(\lambda^{(1)})^2}.$$

All of the diagonal entries in the matrices for $s_2$ are determined by this formula and the values in (6.3).

(c) Let both sides of the equation $s_2^2 = 1$ act on the vector $v_L$ and take the coefficient of $v_L$ in the result. One gets the equation

$$(s_2)_{LM}(s_2)_{ML} + (s_2)_{LL}^2 = 1,$$

where $M$ is the path to $\lambda$ in $\Gamma$ which is not $L$. It follows that

$$(s_2)_{LM}(s_2)_{ML} = (1 + (s_2)_{LL})(1 - (s_2)_{LL}).$$

Because of the freedom in the choice of the off diagonal entries, Proposition (1.12)(b), it follows that we may choose $(s_2)_{ML} = 1 + (s_2)_{LL} = 1 - (s_2)_{MM}$.

∎

**The Iwahori-Hecke algebra** $HG_2(p^2, q^2)$

Let $p, q$ be indeterminates. The Iwahori-Hecke algebra $HG_2(p^2, q^2)$ of type $G_2$ is the associative algebra with 1 over the field $\mathbb{C}(p, q)$ given by generators $T_1, T_2$ and relations

$$
\begin{aligned}
&T_1 T_2 T_1 T_2 T_1 T_2 = T_2 T_1 T_2 T_1 T_2 T_1, \\
&T_1^2 = (p - p^{-1}) T_1 + 1, \quad \text{and} \quad T_2^2 = (q - q^{-1}) T_2 + 1.
\end{aligned}
\tag{6.8}
$$

*Analogues of Jucys-Murphy elements*

If $w_{1,0}$ is the longest element of the Weyl group $S_2 = WA_1$ and $w_{2,0}$ is the longest element of the Weyl group $WG_2$, then the corresponding elements in the Iwahori-Hecke algebras $HA_1(p^2)$ and $HG_2(p^2, q^2)$ are given by

$$
z_1 = T_{w_{1,0}} = T_1, \qquad \text{and} \qquad z_2 = T_{w_{2,0}} = T_1 T_2 T_1 T_2 T_1 T_2.
$$

Following (1.8) and Proposition (2.4), define sets $Z_k = \{z_k\}$, $k = 1, 2$.

*Seminormal representations*

Let $P_1 \subseteq P_2$ be the path algebras over the field $\mathbb{C}(p, q)$ which are associated to the diagram $\Gamma$ in (6.1). For each $\lambda \in \hat{G}_2$, let

$$
V^\lambda = \mathbb{C}(p, q)\text{-}span\{v_L \mid L \in \mathcal{L}(\lambda)\},
\tag{6.9}
$$

so that the vectors $v_L$, indexed by the paths $L = (\square \to \lambda^{(1)} \to \lambda^{(2)} = \lambda)$ in $\Gamma$ which end at $\lambda$, form a seminormal basis of the irreducible $P_2$-module $V^\lambda$. It follows from Lemma (1.9), that for any choice of an isomorphism $\Phi$ between the path algebra $P_2$ and $HG_2(p^2, q^2)$ such that $\Phi(P_1) = \mathbb{C}HA_1(p^2) \subseteq HG_2(p^2, q^2)$,

$$
z_1 v_L = p^{c_{1,0}(\lambda^{(1)})} v_L, \qquad \text{and} \qquad z_2 v_L = c_{2,0}(\lambda^{(2)}) p^{c_{2,s}(\lambda^{(2)})} q^{c_{2,\ell}(\lambda^{(2)})} v_L,
\tag{6.10}
$$

if $L = (\square \to \lambda^{(1)} \to \lambda^{(2)})$ is a path in $\Gamma$.

**Theorem 6.11.** *Irreducible seminormal representations of* $HG_2(p^2, q^2)$

*are given explicitly by*

$$\phi_{1,0}(T_1) = (p), \qquad\qquad \phi_{1,0}(T_2) = (q),$$

$$\phi_{1,6}(T_1) = (-p^{-1}), \qquad\qquad \phi_{1,6}(T_2) = (-q^{-1}),$$

$$\phi'_{1,3}(T_1) = (p), \qquad\qquad \phi'_{1,3}(T_2) = (-q^{-1}),$$

$$\phi''_{1,3}(T_1) = (-p^{-1}), \qquad\qquad \phi''_{1,3}(T_2) = (q),$$

$$\phi_{2,1}(T_1) = \begin{pmatrix} p & 0 \\ 0 & -p^{-1} \end{pmatrix} \quad \phi_{2,1}(T_2) = \begin{pmatrix} a & b \\ c & d \end{pmatrix}$$

$$\phi_{2,2}(T_1) = \begin{pmatrix} p & 0 \\ 0 & -p^{-1} \end{pmatrix} \quad \phi_{2,2}(T_2) = \begin{pmatrix} x & y \\ z & w \end{pmatrix}$$

*where*

$$
\begin{aligned}
a &= \frac{1 + p^{-1}(q - q^{-1})}{p + p^{-1}}, & &\qquad x = \frac{-1 + p^{-1}(q - q^{-1})}{p + p^{-1}}, \\
b &= q - a, & \text{and} &\qquad y = q - x, \\
c &= q^{-1} + a, & &\qquad z = q^{-1} + x, \\
d &= (q - q^{-1}) - a, & &\qquad w = (q - q^{-1}) - x.
\end{aligned}
$$

*Proof.* The proof is entirely similar to the proof of Theorem (6.7). Let us only explain how to get the entries in the matrices $\phi_{2,1}(T_2)$ and $\phi_{2,2}(T_2)$. Let $\lambda = (2,1)$ or $\lambda = (2,2)$ and suppose that $L$ and $M$ are the two paths to $\lambda$ in $\Gamma$. From the definition of the element $z_2$ we get

$$T_1^{-1} T_2^{-1} T_1^{-1} z_2 = T_2 T_1 T_2.$$

By rewriting $T_2^{-1}$ as $T_2 - (q - q^{-1})$ we have

$$T_1^{-1} T_2 T_1^{-1} z_2 - (q - q^{-1}) T_1^{-2} z_2 = T_2 T_1 T_2.$$

Taking the $(L, M)$ entry of each side of the above equation gives

$$(T_1)_{LL}^{-1} (T_2)_{LM} (T_1)_{MM}^{-1} c_{2,0}(\lambda) - 0$$
$$= (T_2)_{LL} (T_1)_{LL} (T_2)_{LM} + (T_2)_{LM} (T_1)_{MM} (T_2)_{MM}.$$

Since these representations are irreducible it follows that $(T_2)_{LM} \neq 0$. Dividing by $(T_2)_{LM}$ and using the fact that $(T_1)_{LL} = p$ and $(T_1)_{MM} = -p^{-1}$, we get the equation

$$-c_{2,0}(\lambda) = p(T_2)_{LL} - p^{-1}(T_2)_{MM}. \tag{i}$$

Now, the equation $T_2^2 = (q - q^{-1})T_2 + 1$ forces that the trace of the matrix of $T_2$ is $q - q^{-1}$, and so

$$(T_2)_{LL} + (T_2)_{MM} = q - q^{-1}. \tag{ii}$$

The values for $x$ and $a$ in the statement of the theorem now follow easily from (i) and (ii). The values of the off diagonal entries are determined (up to a constant, see Proposition (1.12)(a)) by the equation

$$(T_2)_{LM}(T_2)_{ML} = (q^{-1} + (T_2)_{LL})(q - (T_2)_{MM}).$$

This equation is obtained by taking the $(L, L)$ entry in the equation $T_2^2 = (q - q^{-1})T_2 + 1$. ∎